\documentclass[a4paper,12pt]{article}
\usepackage[english]{babel}
\usepackage{amsthm,amsmath,amssymb,amsfonts}
\usepackage{graphicx}
\usepackage{epstopdf}
\usepackage{bbm}

\usepackage[dvips]{color}

\newtheorem {lemma}{Lemma}
\newtheorem {corollary}{Corollary}
\newtheorem {prop}{Proposition}
\newtheorem {thm}{Theorem}

\newcommand{\rmd}{\,\mathrm{d}}

\newcommand{\dist}{\mathop{\mathrm{dist}}}
\newcommand{\NN}{\mbox{${\mathcal N}$}}

\newcommand{\Var}{{\mathrm{Var}}\ }

\def \phi   {\varphi}

\def \a   {\alpha}

\def \g   {\gamma}
\def \d   {\delta}
\def \l   {{\lambda}}
\def \eps {\varepsilon}

\def \ge {\geq}
\def\E{{\mathbb{E}}}

\def\Var{{\mathbb{V}{\sf ar}\,}}
\def\P{{\mathbb{P}}}
\def \R {{\mathbb{R}}}
\newcommand{\Z}{\mathbb{Z}}
\def\F{{\cal{F}}}

\def\limt{\lim_{t\to\infty}}
\def\limt0{\lim_{t\to 0}}

\def\|{\,|\,}
\newcommand{\cd}{\stackrel{\cal D}{\longrightarrow}}
\def\OO{\mathcal{O}}

\def \AA{\mathcal{A}}
\def \I{\mathcal{I}}
\def \II{\mathcal{II}}
\def \p{\partial }

\def\bn#1\en{\begin{align*}#1\end{align*}}
\def\bnn#1\enn{\begin{align}#1\end{align}}
\def\be#1\ee{\begin{equation}#1\end{equation}}
\def\ben#1\een{\begin{equation*}#1\end{equation*}}
\def\beqn#1\eeqn{\begin{eqnarray}#1\end{eqnarray}}
\def\beqnn#1\eeqnn{\begin{eqnarray*}#1\end{eqnarray*}}

\title{Border aggregation model}
\author{Debleena Thacker${}^*$  and Stanislav Volkov\footnote{Centre for Mathematical Sciences, Lund University, Box 118 SE-22100, Lund, Sweden}} 

\begin {document}
\maketitle
\begin {abstract}
Start with a graph with a subset of vertices called {\it the border}. A particle released from the origin  performs a random walk on the graph until it comes to the immediate neighbourhood of the border, at which point it joins this subset thus increasing the border by one point. Then a new particle is released from the origin and the process repeats until the origin becomes a part of the border itself. We are interested in the total number $\xi$ of particles to be released by this final moment. 

We show that this model covers OK Corral model as well as the erosion model, and obtain distributions and bounds for $\xi$ in cases where the graph is star graph, regular tree, and a $d-$dimensional lattice.
\end {abstract}

{{\bf Keywords}: OK Corral model, DLA model, erosion model, random walks, aggregation}

{{\bf Subject classification}:  60K35, 82B24}

\section{Introduction}\label{Intro}
Consider a finite connected graph $G$ (for simplicity $G$ will denote also the set of its vertices) with some designated vertex called {\em origin} $v_0$ and some non-empty set of {\em border vertices} ${\sf B}$. We define recursively set of {\em sticky vertices} $S_n\subseteq G$, with $S_0={\sf B}$.
The process runs as follows: a particle starts some sort of random walk originated at $X_0=v_0$ on $G$, and whenever it comes within the immediate vicinity (i.e.\ one edge away) from a sticky vertex, random walks stops and this particle joins the sets of sticky vertices. Then a new particle starts a random walk  at $v_0$ and runs until it stops, and the process restarts again, until $v_0$ becomes sticky itself, at which point the process stops completely. We are interested in random quantity $\xi=\xi(G)$, the total number of particles emitted from the origin during the lifespan of the process, which always satisfies
$$
\dist(v_0,{\sf B})\le \xi\le |G|-|{\sf B}|
$$ 
where $\dist(A,C)$ for $A,C\subseteq G$ is the number of edges in shortest path connecting $A$ and $C$, and $|G|$ and $|B|$ denote the number of vertices in $G$ and $B$ respectively.

Formally, define a sequence of subsets $S_n$, $n\ge 0$, such that $S_0={\sf B}$ and $S_{n}=S_{n-1}\cup \{w_{n}\}$ where $w_{n}\in G$ is defined as follows: for $n\ge 1$ let $X_t^{(n)}$, $t=0,1,2,\dots$ be a random walk on $G$ such that $X^{(n)}_0=v_0$ and 
$$
\tau_n=\inf\{t\ge 0:\ \dist(X^{(n)}_t,S_{n-1})=1\},
\quad v_n=X^{(n)}_{\tau_n}.
$$
Then
$$
\xi=\inf\{n\ge 1: \dist(v_0,S_{n-1})=1\}\ge 1.
$$
We call the above model border aggregation model ({\sl BA model}  for short).

We study the BA model on a variety of graphs, namely, the star graph, regular $d-$ary trees, and the integer lattice for dimensions $d \geq 2$. Incidentally (Yuval Peres, personal communications; also~\cite{LE&Pe07}), the BA model on a finite piece of the integer line is equivalent to the OK Corral model of~\cite{WM,K99,KV}, where its asymptotic behaviour has been completely analysed. 
Note that BA model was called {\em internal erosion model} in~\cite{LE&Pe07}, however, we feel that the term ``border aggregation model'' is more appropriate. 
The authors of this paper also guess that on a disk of radius $N$ in $\R^2$ the number of {\em eroded points}, which coincides with the number of emitted particles in BA model, grows asymptotically at the rate of $N^{\a}, \mbox{ } \a <2$, the conjecture which we partially solve here.

According to~\cite{LE&Pe07}, the model on $\Z^d$, $d \ge 2$, can be viewed as an ``inversion" of the classical diffusion-limited-aggregation model (DLA), in which particles performing random walks are released at infinity, and they stop once they reach some nearest neighbour of the cluster, which initially consists only of the origin. When $d=2$ Kesten~\cite{Ke87} showed that with probability one the maximum radius of the random cluster is eventually at most of the order $n^{2/3}$, where $n$ is the number of accumulated particles; the corresponding order  for $d\ge 3$ is  $n^{2/d}$. This suggests that for $d=2$ with high probability  the number of emitted particles $\xi$ in the BA model should be at least of the order $N^{3/2}$, and for $d\ge 3$, at least $N^{d/2}$ particles must be emitted. 

In Section~\ref{Sec:2-dim-lattice} we obtain a slightly worse lower bound of $N^{4/3}$ for case $d=2$. In Section~\ref{Sec:d-ge-3} we show that if $d \ge 3$, then $\xi$ must grow at least as $N^{d/2}$ with probability converging to one, which we believe is the correct order.

As it was mentioned above, the BA model is close to the internal diffusion-limited-aggregation-model on $\Z^d$ studied e.g.\ in~\cite{We-Sa-83,LA-Br-Gr-92}. In the process studied in~\cite{LA-Br-Gr-92} initially only the origin is occupied, and particles perform simple  random walks until the moment they visit a vertex outside of the cluster, at which point they stop and become part of the cluster. The authors show that with probability $1$ the limiting shape of the cluster can be well approximated by the $d-$dimensional ball centred at the origin.

These results were further improved in~\cite{Je-Le-Sh-2013, Je-Le-Sh-2014} where it was proved that the maximal error for the DLA cluster is $O( \log t) \text{ and } O( \sqrt{ \log t})$ respectively for $d=2$, and $ d \ge 3$; moreover,  the fluctuations (appropriately normalized) of the cluster converge in law to Gaussian free field. Similar results were also obtained independently in~\cite{As-Ga-2013, As-Ga-2014}. The model has been also studied on the comb lattice in~\cite{As-Ra-2016}, and similar limiting shape theorems have been obtained.  

Another model similar to the BA model has been studied  in~\cite{LaLy99}. The authors introduce a process where one particle is placed at each location in the interval $\left[0, N\right]$, and at every step a randomly chosen particle  from $[1,N]$ is moved to the left until the whole process coalesce at the origin. The authors show that the random time until the coalescence grows asymptotically at the rate of $N^{3/2}$, and the variance is upper bounded by $N^{5/2}$.

The rest of the paper is organized as follows. In Section~\ref{Sec:star} we study the BA model on star graphs. In Section~\ref{Sec:trees} we obtain quite sharp bounds for $\xi$ on regular trees. 
In Section~\ref{Sec:2-dim-lattice} we analyze the model on the two-dimensional lattice and a comb lattice; in Section~\ref{Sec:d-ge-3} we obtain the results for higher dimensions. In the latter two cases we obtain  non-trivial lower bounds only. 

Finally we mention that throughout the paper for any two positive sequences $\{a_N\}_{N \ge 1}$ and $\{b_N\}_{N \ge 1}$ $a_N \sim b_N$ denotes the fact  that $\lim_{N \to \infty}{a_N}/{b_N}=1.$

\section{Star graph}\label{Sec:star}
\setlength{\unitlength}{1cm}
\begin{picture}(16,4)
\put(7.8,2.3){$v_0$}
\put(8,2){\line(3,1){5}}
\put(8,2){\line(3,-1){5}}
\put(8,2){\line(1,0){5.5}}
\put(8,2){\line(-4,-1){5}}
\put(8,2){\line(-4,1){5}}
\multiput(8, 2)(0.6,  0.2 ){9}{\circle*{0.15}}
\multiput(8, 2)(0.6, -0.2 ){9}{\circle*{0.15}}
\multiput(8, 2)(0.64, 0   ){9}{\circle*{0.15}}
\multiput(8, 2)(-0.6,-0.15){9}{\circle*{0.15}}
\multiput(8, 2)(-0.6, 0.15){9}{\circle*{0.15}}
\put(13,2){\oval(0.5,4)}\put(13.5,3){${\sf B}$}
\put(3.1,2){\oval(0.5,4)}\put(2.3,3){${\sf B}$}
\put(7,1){$K=5$}
\end{picture}

Let $G$ consist of $K\ge 2$ pieces of $\{0,1,2,\dots,N+1\}\subset \Z_+$ sharing a common origin $v_0=0$, and let ${\sf B}$ be the $K$ endpoints $(N+1)$ of each of the segment. Let $X$ be a simple random walk on $G$. 

If $K=2$ then $G=\{-(N+1),-N,\dots,-1,0,1,\dots,N,N+1\}$ with $v_0=0$ and ${\sf B}=\{-N-1\}\cup\{N+1\}$. Let $X$ be a simple random walk on $G$. As it was mentioned in the Introduction, this model is equivalent to the OK Corral model of~\cite{WM} with the initial number of shooters equal to $N$, so $2N-\xi$ gives the number of survivals in the positive or negative group, and it is asymptotic  order is $N^{3/4}$, as it was shown in~\cite{KV}, where the distribution was found explicitly.
The case $K\ge 3$ is thus a natural generalization of the OK Corral model.

Using the elementary properties of simple random walk it is easy to deduce that the border aggregation model on this star graph is equivalent to the following urn-like model. Let $X_i(0)=N$, $i=1,2,\dots,K$.  Given the vector $X(j)=(X_1(j),\dots,X_K(j))$, $j=0,1,2,\dots$, we independently sample $\zeta(j)\in\{1,2,\dots,K\}$ such that
$$
\P(\zeta(j)=k)=\frac{X_k(j)^{-1}}{\sum_{i=1}^K X_i(j)^{-1}}
$$
and let
$$
X(j+1)=X(j)-(0,0,\dots,0,\underbrace{1}_{\zeta(j)\text{ place}},0,\dots,0).
$$
In other words, at each moment of time exactly one of the $X_i$'s decreases by $1$, and the chances of the $k$-th segment to be picked are inversely proportional to its length. Again, in case $K=2$ this is exactly the OK Corral model, and we can think of this model as a generalization of the latter one with $K$ different groups. Let $S_N(K)$ be the number of survivals by the time one of the group is eliminated; then $S_N(K)=NK-\xi$. We will show that $S_N(K)\sim N^{3/4}$ for any fixed value of $K$. 

The crucial observation here is that we can couple the process with $K$ independent continuous time processes in the same fashion as it is done in~\cite{KV} using the idea of Rubin's construction from~\cite{BD}. Indeed, let $Y_{i}(t)$, $i=1,2,\dots,K$, $t\ge 0$, be $K$ independent pure death process all starting at $Y_i(0)=N$ with the death rate at level $k$ equal to $1/k$. Fix $K \ge 2$ throughout. Let 
\bnn\label{eq:istar}
\tau_{i}&=\inf\{t:\ Y_{i}(t)=0\},\quad 
\bar\tau =\min_{i=1,2,\dots,K} \tau_{i}=\tau_{i^*},\quad
i^*=\arg\min\tau_{i} \nonumber \\
S_N(K)&=\sum_{i=1}^K Y_{i}(\bar\tau)=
\sum_{i=1,i\ne i^*}^K Y_{i}(\bar\tau)
\enn
i.e., $i^*$ is the index of the process which dies out first. It is clear from the context that $\bar\tau$ depends on $K,$ but we simply write $\bar\tau$, instead of $\bar\tau(K),$ as this does not create any ambiguity. Then the definition of $S_N(K)$ is consistent with the definition given earlier in this section. 

Observe that for any $i \geq 1$, there exist independent random variables $\left\{\xi_{i,k}\, , k=1, 2, N\right\}$ such that $\xi_{i,k}\sim Exp\left(1/k\right),$ such that 
\be\label{eq:tauexp}
\tau_i= \displaystyle \sum_{k=1}^{N} \xi_{i,k}.
\ee 

Then $\tau_i$ satisfies the following Central Limit Theorem (CLT). 
\begin{thm}\label{thm:CLTtaui} Let $\tau_i$ be as defined above. Then, as $N \to \infty$
\be\label{eq:CLTtaui}
\frac{\tau_i-\frac{N^2}{2}}{\frac{N^{\frac{3}{2}}}{\sqrt{3}}}\cd \NN\left(0,1\right)
\ee 
where $\cd$ denotes convergence in distribution.
\end{thm}
\begin{proof} It is easy to see from~\eqref{eq:tauexp} 
\bn
\E\left[\tau_i\right]&= \sum_{k=1}^{N}k=\frac{N (N+1)}{2}\sim \frac{N^2}2 ,
\\
\Var\left(\tau_i\right)&= \sum_{k=1}^N k^2=\frac{N(N+1)(2N+1)}{6} \sim \frac{N^3}3.
\en
The rest of the proof now follows from an easy application of the standard CLT for the sum of independent random variables.
\end{proof}

\begin{thm}\label{Theorem:B-E_bound_tau_i}
Let $\tau_i$ be as defined above. Then, as $N \to \infty$
\ben\label{Eq:B-E_bound_tau_i}
 \sup_{x \in \R}  \left\vert 
\P \left(\frac{\tau_i-\frac{N^2}{2}}{\frac{N^{\frac{3}{2}}}{\sqrt{3}}} \le x\right)
-\Phi(x)\right\vert\le (2.75) l_{3,N}=\OO\left(\frac{1}{\sqrt N}\right),
\een 
where  
\ben\label{Eq:def_l_3}
l_{3,N}:= \frac{\frac{1}{N}\sum_{j=1}^N \E \big{\vert} \xi_{1,j}-j\big{\vert}^3 }{\left(\frac{1}{N}\sum_{j=1}^N \Var(\xi_{1,j})\right)^{3/2}}\, N^{-1/2}.
\een

\end{thm}
\begin{proof}
The proof is an easy application of Theorem~12.4 (Berry-Essen Theorem) on p.~104 in~\cite{BR}.
\end{proof}
\begin{thm}\label{Thm:Local_CLTtaui}
Let $f_{N,i}(x)$ 
denote the density function of 
$$
\frac{\tau_i-N(N+1)/2}{N^{3/2}/\sqrt 3}=
\tau_i \sqrt{\frac{3}{N^3}}- z_N
$$ 
where 
$z_N=\frac{(N+1)\sqrt 3}{2 \sqrt{N}}
=\sqrt{3N} \left(\frac{1}{2}+\frac{1}{2N}\right)
$
and $n(x)=\frac{\exp(-x^2/2)}{\sqrt{2\pi}}$ is the density function of $\NN(0,1)$. Then for every $\epsilon>0$ 
\ben\label{Eq:Local_CLTtaui}
\sup _{x \in \R} \vert f_{N,i}(x) - n(x) \vert = o \left(N^{-(1/2-\epsilon)}\right)
\een
as $N\to\infty$.
\end{thm}
The proof of this theorem is similar to the proof of Lemma~2 of~\cite{KV}.
\begin{proof}
By the Fourier inversion formula (e.g.\ Theorem 3.3.5 in~\cite{DUR}), 
\be \label{Eq:Fourier_inversion}
f_{N,i}(x) - n(x) =\frac{1}{2\pi}\int_{-\infty}^{\infty}{e^{-itx}\left(\phi_{N}(t)-e^{-t^2/2}\right)\, \rmd t}
\ee 
where 
\bn
 \phi_N(t)  :=  \E e^{it\, \frac{\tau_i-{N(N+1)}/{2}}{{N^{3/2}}/{\sqrt 3}}}
=  \E e^{it\tau_i \, \sqrt{\frac{3}{N^3}}-it z_N}
  = e^{-itz_N} 
 \prod _{k=1}^N 
 \left[1-\frac{ikt}{{N^{3/2}}/{\sqrt 3}}\right]^{-1}
\en
is the characteristic function for $\frac{\tau_i-N(N+1)/2}{N^{3/2}/\sqrt{3}}$, satisfying $\int |\phi_N(t)|\rmd t<\infty$ once $N\ge 2$.

For $\vert t \vert \ll \sqrt{N}$, using series expansion $-\ln(1-\a)=\a+\a^2/2+\a^3/3+\dots$ we have 
\bn
 itz_N +\ln \phi_N (t) & =  - \sum_{k=1}^N \ln \left(1-\frac{ikt\sqrt{3}}{N^{3/2}}\right)  =   \sum_{k=1}^N \sum_{j=1}^{\infty} \frac{(itk)^j   3^{j/2}}{j N^{3j/2} } =   \sum_{j=1}^{\infty}  
A_{N,j}  \frac{(it)^j 3^{j/2} } {j N^{j/2-1 }} 
\en 
where 
$$
A_{N,j}=
\frac{1}{N^{j+1}}\sum_{k=1}^N k^j=
\begin{cases}
\frac 12 +\frac 1{2N},& j=1,\\
\frac 1{j+1} +\frac 1{2N}+\OO(N^{-2}),& j\ge 2.
\end{cases}
$$
by interchanging the order of summation. Hence 
\bn
\ln \phi_N(t)
&= 
- itz_N+\left(\frac{1}{2}+\frac{1}{2N}\right)
it\sqrt{3N}+
 \sum_{j=2}^{\infty}  
\left[\frac 1{j+1}+\frac{1}{2N}+\OO(N^{-2})\right]
 \frac{(it)^{j} 3^{j/2} } {j N^{j/2-1 }} 
\\ 
&
= \left[ -\frac{t^2}{2} 
- \frac{i\,\sqrt{3}}{4}\frac{t^3}{N^{1/2}}
+ \frac{9}{20}\, \frac{t^4}{N}
+ \frac{3i\sqrt{3}}{10}\, \frac{t^5}{N^{3/2}}
+\dots \right]
+\OO(t^2/N)
\\
&=
-\frac{t^2}{2}
\left\{1
+ \frac{i\,\sqrt{3}}{2} \left(\frac{t}{\sqrt N}\right)
- \frac{9}{10}\, \left(\frac{t}{\sqrt N}\right)^2
+\dots \right\}+\OO(t^2/N).
\en

As in the proof of Lemma~2 in~\cite{KV},  we divide the area of integral in \eqref{Eq:Fourier_inversion} into two parts; $[-N^{\delta}, N^{\delta}],$ where $0< \delta < \min (1/2, \epsilon/4)$ and its complement. 
\bnn \label{Eq:Local_CLTtau_i_logN}
& \left|\int_{-N^{\d}}^{N^{\d}}{e^{-itx}(\phi_{N}(t)-e^{-\frac{t^2}2}) \rmd t}\right| 
  \nonumber
  \le \int_{-N^{\d}}^{N^{\d}}\left| e^{-itx-\frac{t^2}2}\right|\cdot \left| e^{-\frac{t^2}{2}\left[ \frac{i\,\sqrt{3}}{4} \left(\frac{t}{\sqrt N}\right) (1+o(1))\right]}-1\right| \rmd t
\\
& \le \int_{-N^{\d}}^{N^{\d}}
\left|  \frac{t^3\,\sqrt{3}}{8\sqrt N}
 (1+o(1)) \right| \rmd t
= \OO\left(N ^{-(1/2- 4 \delta)}\right).
\enn
The integral over the complement is also dealt in similar fashion as in Lemma~2 in~\cite{KV}. 
\begin{align}\label{Eq:Local_CLTtau_i_other_part} 
\left|\int_{ |t|\ge N^\delta }{e^{-itx}(\phi_{N}(t)-e^{-\frac{t^2}2})\, \rmd t} \right| & \le &  
\int_{  |t|\ge N^\delta  }{(\vert \phi_{N}(t) \vert \, +e^{-\frac{t^2}2 }) \, \rmd t}= o \left(N^{-1}\right),
\end{align}
where the last equality can be obtained by using bounds similar to~equations~(12) and~(14) of~\cite{KV}, using the fact that 
$$
|\phi_{N}(t)|^2 = \prod_{k=1}^N 
\frac 1{1+3 k^2 t^2/N^3}
=\frac{1}{1+t^2[1+\OO(N^{-1})]+
t^4[1/2+\OO(N^{-1})]+\dots}.
$$ 
Now from~\eqref{Eq:Local_CLTtau_i_logN} and~\eqref{Eq:Local_CLTtau_i_other_part} it follows that for all $x \in \R,$
\begin{align*}
\vert f_{N,i}(x) - n(x) \vert = \OO\left(N ^{-(1/2- 4 \delta)}\right) + o \left(N^{-1}\right)=o \left(N^{-(1/2-\epsilon)}\right).
\end{align*}
\end{proof}  

The CLT  and Theorem~\ref{Thm:Local_CLTtaui} implies that  $\tau_i$ has the following representation
$$
\tau_i=\frac{N^2}{2}+\frac{N^{3/2}}{\sqrt{3}} \eta_i+ \OO(N)
$$
where $\eta_i$ are i.i.d.\  normal $N(0,1)$ random variables.

Our first observation is that for each $K$ there exists constants $C(\nu, K)>0,$ such that 
\ben\label{Eq:mu_nu_bdd}
\displaystyle \limsup_{N \to \infty}\frac{
\E \left[ \left(S_N(K)\right)^\nu\right]
}
{(N^{3/4})^{\nu}} 
\le C(\nu, K).
\een  

Indeed, for $K=2$ we already know the result from~\cite{K99}
\be\label{Eq:mu_nu_2}
\displaystyle \lim_{N \to \infty}\frac{\E \left[ \left(S_N(2)\right)^\nu\right]}{(N^{3/4})^{\nu}}= \left(2/3\right)^{\nu}\frac{\Gamma(\frac{\nu}{4}+\frac{1}{2})}{\Gamma(1/2)}=
\left(2/3\right)^{\nu}\frac{\Gamma(\frac{\nu}{4}+\frac{1}{2})}{\sqrt{\pi}}=:C(\nu,2).
\ee
We will use this limit to obtain the result for general $K \geq 3$. Note that by symmetry $i^*$ has a uniform distribution over $\{1,\dots,K\}.$

Let $\a:=\left(\a_1, \a_2, \dots, \a_{K-1}\right)\in\Z_+^{K-1}$, $|\a|=\a_1+\dots+\a_{K-1}$, and ${\nu\choose\a}=\frac{\nu!}{\a_1!\a_2!\dots\a_{K-1}!}$ be the multinomial coefficient. Then
\begin{align*}
\E \left[\left(S_N(K)\right)^{\nu}\right] 
& =  \E \left[
\sum _{i^*=1}^K 
\left(\sum _{i=1, i\ne i^*}^K Y_i(\bar\tau)\right)^{\nu}
\right]
 =  K\, \E \left[\left(\sum _{i=1, i^*=K}^{K-1} Y_i(\bar\tau)\right)^{\nu}\right]\\
& =
 \sum_{\a:\ |\a|=\nu} {\nu \choose \a} \, \E \left[\prod _{i=1}^{K-1} Y_i^{\a_i}(\bar\tau) \mathbf{1}_{\{\bar\tau= \tau_K\}}\right]
\\
&
 = \sum_{\a:\ |\a|=\nu} {\nu \choose \a} \prod _{i=1}^{K-1}\E \left[Y_i^{\a_i}(\bar\tau) \mathbf{1}_{\{\bar\tau= \tau_K\}}\right]
\end{align*}
using the symmetry, and the independence of each of $Y_i(\bar\tau)$ conditioned on the event $\{\bar\tau= \tau_K\}.$ It is easy to see that for any $i$
\bn
(*)=\E \left[Y_1^{\a_i}(\bar\tau) \mathbf{1}_{\{\bar\tau= \tau_K\}}\right] 
& = \int_0^{\infty}{Y_1^{\a_i}(s)\P\left(\tau_1>s\right)\dots \P\left(\tau_{K-1}>s\right)f_{\tau_K}(s) \, ds}\\
& \le  \int_0^{\infty}{Y_1^{\a_i}(s)\P\left(\tau_1>s\right)\, f_{\tau_K}(s) \, ds}
\en
where $f_{\tau_K}$ denotes the density of $\tau_K$.
Now note that $\bar\tau(K)$ stochastically decreases in $K$, therefore the expression on the RHS is smallest when $K=2$, therefore, assuming $K=2$ in the RHS, we get
\bn
(*) \le  \int_0^{\infty}{Y_1^{\a_i}(s)\, \P\left(\tau_1>s\right) f_{\tau_2}(s) \, ds}
=\frac{1}{2} \E \left[\left(S_N(2)\right)^{\a_i}\right],
\en
where the last equality follows from symmetry between $Y_1$ and $Y_2$. Consequently, using the limit in~\eqref{Eq:mu_nu_2}, we have 
\bn
\displaystyle 
\limsup_{N \to \infty} \frac{\E \left[ \left(S_N(K)\right)^\nu\right]}{(N^{3/4})^{\nu}} \le  \left(\frac23 \right)^{\nu/4}
\pi^{-\frac{K-1}{2}}
\sum_{\a:\ |\a|=\nu} {\nu \choose \a} \prod _{i=1}^{K-1} \Gamma\left(\frac{\a_i}{4}+\frac{1}{2}\right).
\en

\begin{corollary}\label{cor.unif.bound}
The sequence of $\displaystyle\frac{S_N(K)}{N^{3/4}}$ is uniformly bounded in ${\mathbb L}^\nu$ for every $\nu=1,2,\dots$.
\end{corollary}

Introduce integer-valued function $m_j: \{1,2,\dots, K-1\}\to \{1,2, \dots, K\}\setminus\{j\} $ such that
\bn
m_j(i) = \begin{cases} i, \mbox{ if } i <j, \\
                   i+1, \mbox{ if } i >j,
\end{cases} 
\en
and let
\bn
\zeta_i^{(N)}:=\frac{Y_{m_{i^*}(i)}(\bar\tau)}{N^{3/4}}>0,
\quad i=1,2,\dots,K-1,
\en
where $i^*$ is defined in~\eqref{eq:istar}, be the lengths of the $K-1$ rays which have not been filled in by the time $\bar\tau$.

\newpage
\begin{thm}\label{Theorem:convergence_joint_distribution_K_rows}
We have
\ben\label{Eq:convergence_joint_distribution_K_rows}
 \zeta^{(N)} \cd \zeta=\left(\zeta_1,\zeta_2,\dots,\zeta_{(K-1)} \right),
\een where $\zeta$ is a non-degenerate jointly continuous random variable satisfying
\bn
\P \left(\zeta_1>a_1,\dots,\zeta_{K-1}> a_{K-1} \right)
=\frac{K}{\sqrt{2 \pi}}\int_{-\infty}^{\infty}{\prod_{i=1}^{K-1}\left[1 -\Phi\left(\frac{\sqrt{3}}{2}a_i^2+w\right)\right]\,e^{-\frac{w^2}{2}}\, \rmd w}
=:G(\mathbf{a})
\en 
for any  $\mathbf{a}=(a_1,a_2,\dots,a_{K-1})\in \R^{K-1}_+$. Moreover, the joint density of $\zeta$ is given by
$$
f_{\zeta}(\mathbf{a})=a_1 a_2\dots a_{K-1}\,
\sqrt{ K \left[\frac 3{2 \pi}\right]^{K-1} }
\exp\left\{-
\frac{3}{8}\left[ \sum_{i=1}^{K-1} a_i^4-\frac1K \left(\sum_{i=1}^{K-1} a_i^2\right)^2
\right]
\right\}.
$$
Therefore $\frac {\sqrt 3}2\, \zeta^2=\left(\frac {\sqrt 3}2\, \zeta_1^2,\dots,\frac {\sqrt 3}2\, \zeta_{K-1}^2\right)$ has the density
$$
f(\mathbf{b})=
\sqrt{ 
\frac K{\left(2\pi\right)^{K-1}} 
}\,
\exp\left\{-
\frac{1}{2}\left[ \sum_{i=1}^{K-1} b_i^2-\frac1K \left(\sum_{i=1}^{K-1} b_i\right)^2
\right]
\right\}.
$$
for $\mathbf{b}\in \R^{K-1}_+$ and
thus we can write $\frac {\sqrt 3}2\, \zeta^2$ as 
$\sqrt{\frac K{K-1}}\cdot Z$ conditioned on $Z_1\ge 0,\dots,Z_{K-1}\ge 0$
where $Z_i=W_i-\frac {\sum_{i=1}^{K-1} W_i}{\sqrt K -1}$ 
and $\{W_i\}_{i=1}^{K-1}$ are iid $N(0,1)$ (thus
 ${\sf Cov}(Z_i,Z_j)=1+\mathbf{1}_{\{i=j\}}$, $i,j=1,\dots,K-1$).
\end{thm}

\begin{corollary}\label{Cor:convergence_single_component}
As $N\to\infty$ we have:
\begin{itemize}
\item[(a)]
$\displaystyle \frac{Y_{i}(\bar\tau)}{N^{3/4}}\cd \eta$ 
where  $i=1,2,\dots,K$ and the CDF $F_\eta(x)$  is
\bn
\begin{cases}
1- \frac{K-1}{\sqrt{2 \pi}}\int_{-\infty}^{\infty}
\left[1 -\Phi\left(\frac{\sqrt{3}}{2}x^2+w\right)\right] \left[1-\Phi(w)\right]^{K-1} e^{-\frac{w^2}{2}}\, \rmd w, &\text{ if $x>0$;}\\
1/K,  &\text{ if $x=0$;}\\
0,  &\text{ if $x<0$,}\\
\end{cases}
\en
that is, $\eta$ is a mixture of an atom at $0$ and a continuous distribution on $\R_+$;
\item[(b)] As $N\to\infty$ we have
$\displaystyle \frac{S_N(K)}{N^{3/4}} \cd \sum_{i=1}^{K-1} \zeta_i.$ and morevoer $\displaystyle \E\left[\frac{S_N(K)}{N^{3/4}} \right]^\nu\to  \E\left[\sum_{i=1}^{K-1} \zeta_i\right]^\nu$ for any positive integer $\nu$.
\end{itemize}
\end{corollary}
\begin{proof}
Part (a) immediately follows from Theorem~\ref{Eq:convergence_joint_distribution_K_rows} and the definition of~$\zeta^{(N)}$. To show part (b), define the function $g \colon \R^{K-1} \to \R$ as
$g(x):=\sum_{i=1}^{K-1}x_i$ for $x=(x_1,\dots,x_{K-1})$.
It is easy to see that $g(\cdot)$ is continuous and 
$g\left(\zeta^{(N)}\right)={S_N(K)}/{N^{3/4}}.$
From~\eqref{Eq:convergence_joint_distribution_K_rows} and the continuous mapping theorem, see~\cite{DUR} Theorem~3.2.4, we have  
$
g\left(\zeta^{(N)}\right)\cd g(\zeta)=\sum_{i=1}^{K-1} \zeta_i.
$

Finally, the statement about convergence of expectations follows from our Corollary~\ref{cor.unif.bound} and  Corollary on p.~348 of~\cite{BIL}.
\end{proof}

\begin{proof}[Proof of Theorem \ref{Theorem:convergence_joint_distribution_K_rows}]
By symmetry
\begin{align*}
  \P\left(\zeta_1>a_1,\dots,\zeta_{K-1}> a_{K-1} \right) 
& =  \sum_{j=1}^K \P\left(\zeta_1>a_1,\dots,\zeta_{K-1}> a_{K-1}, \, i^*=j \right)
 \\  
 & =   \sum_{j=1}^K \P\left(\frac{Y_i(\bar\tau)}{N^{3/4}}>a_{m_{j}^{-1} (i)}\, \forall i \neq j; \, i^*=j  \right)
\\ 
&=  K\, \P\left(\frac{Y_1(\bar\tau)}{N^{3/4}}>a_1,  \dots, \frac{Y_{K-1}(\bar\tau)}{N^{3/4}}>a_{K-1}, \, i^*=K  \right).
\end{align*}
Define the stopping times, $\tau_{i}=\tau_i(a):= \inf\{t:\ Y_{i}(t)\le a_i N^{3/4}\}$ for  $i =1,2,\dots, K$, and recall that $\{i^*=K\}=\{\bar\tau = \tau_K \}.$ It is easy to see that
\beqn \label{Eq:zeta_y_tau_integrals}
\nonumber & & \P\left(\frac{Y_1(\bar\tau)}{N^{3/4}}>a_1, \frac{Y_2(\bar\tau)}{N^{3/4}}>a_2, \dots, \frac{Y_{(K-1)}(\bar\tau)}{N^{3/4}}>a_{K-1}, \, i^*=K \right)\\
& = &\int_{-\infty}^{\infty}{\P\left(\tau_1>s\right)\, \P\left(\tau_2>s\right)\dots \P\left(\tau_{K-1}>s\right)f_{\tau_K}(s) \, ds}
\\ \nonumber
&=& \I_N+\II_N
\eeqn 
where $f_{\tau_K}(\cdot)$ is the density of the stopping time $\tau_K$,
\bn
\I_{N}&:=\int_{\AA_{N}}{\P\left(\tau_1>s\right)\, \P\left(\tau_2>s\right)\dots \P\left(\tau_{K-1}>s\right)f_{\tau_K}(s) \, \rmd s},
\\
\II_N&:=\int_{\AA^c_{N}}{\P\left(\tau_1>s\right)\, \P\left(\tau_2>s\right)\dots \P\left(\tau_{K-1}>s\right)f_{\tau_K}(s) \, \rmd s},
\\
\text{and } \AA_{N}&:= \left[\frac{N^2}{2}-  \frac{N^{3/2}\ln N}{\sqrt{3}}, \, \frac{N^2}{2}+ \frac{N^{3/2}\ln N}{\sqrt{3}}\right] .
\en
We will now estimate the integral $\I_N.$
Using the representation of $Y$ with the exponential random variables $\xi_{i,k}$ we have
\bn 
\E[\tau_{i}] & = \sum_{k=\lfloor a_i N^{3/4}\rfloor}^N \E [\xi_{i,k}] \sim \frac{N(N+1)-a_i^2 N^{3/2}}{2}; 
\\
\Var\left(\tau_i\right) &=  \sum_{k=\lfloor a_i N^{3/4}\rfloor}^N
 \Var( \xi_{i,k} ) \sim \frac{N^3-a_i^3 N^{2.25}}{3}
 \sim \frac{N^3}{3}.
\en 
Similar to \eqref{eq:CLTtaui}, it is easy to see that 
\be \label{Eq:CLT_tau_a_i}
\frac{\tau_{i,a_{N,i}}-\frac{N(N+1)-a_i ^2 N^{3/2}}2}
{\frac{N^{3/2}}{\sqrt{3}}}\cd \NN (0,1).
\ee 
Applying the change of variables $s=\frac{N(N+1)}2+w\frac{N^{3/2}}{\sqrt 3 }$   for $w \in [-t, t]$ in $\I_N$ we get
\bn
 \I_N= 
 \frac{N^{3/2}}{\sqrt 3 }
   \int_{-\ln N}^{\ln N}
   &  \left[\prod_{j=1}^{K-1}\P\left(\tau_{j,a_{N,j}}>\frac{N(N+1)}{2}+ w\, \frac{N^{3/2}}{\sqrt 3 }\right)\right]
\\ &
\times f_{\tau_K}\left(\frac{N(N+1)}{2}+w\, \frac{N^{3/2}}{\sqrt 3 }\right) \rmd w.
\en
From~\eqref{Eq:CLT_tau_a_i}, it follows that 
\bn 
\P\left(\tau_j>\frac{N(N+1)}{2}+w \frac{N^{3/2}}{\sqrt{3}}\right)&=\P \left(\frac{\tau_j-\frac{N(N+1)-a_j ^2 N^{3/2}}2}{\frac{N^{3/2}}{\sqrt 3}}> \frac{\sqrt 3}2 a_j^2+w \right)
\\ &
 \to 1 -\Phi \left(\frac{\sqrt 3}2 a_j^2+w\right).
\en 
Fix a small $\eps>0$. From Theorem~\ref{Thm:Local_CLTtaui}, we have 
\ben
\sup_{w \in \R} \left|
 \frac{N^{3/2}}{\sqrt 3}\, f_{\tau_K}
\left(\frac{N(N+1)}2+w \frac{N^{3/2}}{\sqrt 3}\right)-n(w) \right|= o \left(N^{-(1/2-\epsilon)}\right).
\een
By dominated convergence theorem, we have for all $N$ large enough 
\bnn \label{Eq:I}
\left|\I_N-\frac{1}{\sqrt{2 \pi}}\int_{-\ln N}^{\ln N} {\prod_{i=1}^{K-1}\left[1 -\Phi\left(\frac{\sqrt{3}}{2}a_i^2+w\right)\right]\,e^{-\frac{w^2}{2}}\, d w}\right| = o\left(N^{-(1/2-\epsilon)}\ln N\right)
\enn
and at the same time
\bnn \label{Eq:Ibis}
& 
\left|
 G(\mathbf{a})
 -\frac 1{\sqrt{2 \pi}}\int_{-\ln N}^{\ln N} {\prod_{i=1}^{K-1}\left[1 -\Phi\left(\frac{\sqrt{3}}{2}a_i^2+w\right)\right]\,e^{-\frac{w^2}{2}}\, d w}\right| 
 \\ \nonumber & \qquad \qquad  \le  
\int_{\ln N}^{\infty} e^{-\frac{w^2}2}\rmd w= o\left(N^{-1}\right).
\enn
By Theorem~\ref{Theorem:B-E_bound_tau_i}
\bnn\label{Eq:II}
|\II_N| & \le \int_{\AA^c_N} {f_{\tau_K}(s) \, ds}=\P \left(\tau_K\ge \frac{N^2}{2}+ \frac{N^{\frac 32}\ln N}{\sqrt 3}\right)+ \P \left(\tau_K \le \frac{N^2}{2}- \frac{N^{\frac 32}\ln N}{\sqrt{3}}\right) \nonumber
\\ &
= 2 (1-\Phi(\ln N))+ \OO\left(\frac{1}{\sqrt{N}}\right)=\OO\left(\frac{1}{\sqrt{N}}\right). 
\enn
Combining~\eqref{Eq:I}, \eqref{Eq:Ibis} and~\eqref{Eq:II} we get the desired convergence in distribution.

Finally, we need to show that the limiting random variable $\zeta$ is jointly continuous and find its density. Since all the partial derivatives of the expression inside the integral sign in the definition of $G(\mathbf{a})$ are continuous, we can interchange integration and differentiation by Leibniz integral rule to obtain that $\zeta$ has the joint density
\begin{align*}
& (-1)^{K-1}\,\frac{\p^{K-1} G(\mathbf{a})}{\p a_1\dots\p a_{K-1} } =\frac{(-1)^{K-1} K}{\sqrt{2 \pi}}\int_{-\infty}^{\infty}{\prod_{i=1}^{K-1}\frac{\p }{\p a_i}\left[1 -\Phi\left(\frac{\sqrt{3}}{2}a_i^2+w\right)\right]\,e^{-\frac{w^2}{2}}\, \rmd w}
\nonumber \\
& \qquad =
\frac{a_1 a_2\dots a_{K-1}\,K\, 3^{\frac{K-1}{2}}}{(2 \pi)^{K/2} }\int_{-\infty}^{\infty}
\exp\left\{-\frac{K w^2+w\sqrt{3}\sum_{i}a_i^2+\frac34 \sum_i a_i^4}
2
\right\}\, \rmd w
\nonumber \\
&  \qquad =
\frac{a_1 a_2\dots a_{K-1}\,\sqrt{K} }{(2 \pi/3)^{(K-1)/2} }
\exp\left\{-
\frac{3}{8}\left[ \sum_i a_i^4-\frac1K \left(\sum_i a_i^2\right)^2
\right]
\right\}
\end{align*}
where the sum is taken over $i=1,2,\dots,K-1$.
\end{proof}

\section{Binary tree (and other regular trees)}
\label{Sec:trees}
\setlength{\unitlength}{8mm}
\begin{picture}(16,5)
\put(8,0.3){$v_0$}
\put(8,0){\vector(4,1){4}}
\put(8,0){\vector(-4,1){4}}
\put(4,1.05){\vector(2,1){2}}
\put(4,1.05){\vector(-2,1){2}}
\put(12,1.05){\vector(2,1){2}}
\put(12,1.05){\vector(-2,1){2}}
\put(2,2.1){\vector(1,1){1}}
\put(2,2.1){\vector(-1,1){1}}
\put(6,2.1){\vector(1,1){1}}
\put(6,2.1){\vector(-1,1){1}}
\put(10,2.1){\vector(1,1){1}}
\put(10,2.1){\vector(-1,1){1}}
\put(14,2.1){\vector(1,1){1}}
\put(14,2.1){\vector(-1,1){1}}
\put(1, 3.15){\vector(1,2){0.5}}\put(1, 3.15){\vector(-1,2){0.5}}
\put(3, 3.15){\vector(1,2){0.5}}\put(3, 3.15){\vector(-1,2){0.5}}
\put(5, 3.15){\vector(1,2){0.5}}\put(5, 3.15){\vector(-1,2){0.5}}
\put(7, 3.15){\vector(1,2){0.5}}\put(7, 3.15){\vector(-1,2){0.5}}
\put(9, 3.15){\vector(1,2){0.5}}\put(9, 3.15){\vector(-1,2){0.5}}
\put(11, 3.15){\vector(1,2){0.5}}\put(11, 3.15){\vector(-1,2){0.5}}
\put(13, 3.15){\vector(1,2){0.5}}\put(13, 3.15){\vector(-1,2){0.5}}
\put(15, 3.15){\vector(1,2){0.5}}\put(15, 3.15){\vector(-1,2){0.5}}
\put(0.5,4.1){\circle{0.3}}\put(1.5,4.1){\circle{0.3}}
\put(2.5,4.1){\circle{0.3}}\put(3.5,4.1){\circle{0.3}}
\put(4.5,4.1){\circle{0.3}}\put(5.5,4.1){\circle{0.3}}
\put(6.5,4.1){\circle{0.3}}\put(7.5,4.1){\circle{0.3}}
\put(8.5,4.1){\circle{0.3}}\put(9.5,4.1){\circle{0.3}}
\put(10.5,4.1){\circle{0.3}}\put(11.5,4.1){\circle{0.3}}
\put(12.5,4.1){\circle{0.3}}\put(13.5,4.1){\circle{0.3}}
\put(14.5,4.1){\circle{0.3}}\put(15.5,4.1){\circle{0.3}}
\put(8,4.1){\oval(16,0.6)}\put(8,4.6){${\sf B}$}\put(8,2){$K=4$}
\end{picture}

Let $G=G_{d,K}$ be a regular $d$-ary tree ($d\ge 2$) with root $v_0$ truncated at level $K$, that is it is $v_0$ and all the vertices of distance no more than $K$ from the origin; thus $|G|=1+d+d^2+\dots+d^K$. We assume that all $d^K$ most remote vertices are the border, and the random walk moves only upwards (away from the root) with equal probability.

Let us now assume that $d=2$ and for the rest of this section  deal only with the binary rooted tree, unless said otherwise. Let $\xi_{K}$ denote the total number of emitted particles on $G_{2,K}$ until $v_0$ becomes a part of the border. Then
$$
\xi_{1}=1,\quad  \xi_2=2, \quad
\xi_3=\begin{cases}
3&\text{with probability }1/2,\\
4&\text{with probability }1/2
\end{cases}
$$
and in general
\bnn\label{eqrec}
\xi_{K+1}=1+\eta(\xi_K',\xi_K'')
\enn
where $\xi_K'$ and $\xi_K''$ are two independent copies of $\xi_K$ and $\eta(a,b)$ is the number of tosses of a fair coin required to reach either $a$ heads or $b$ tails, whichever comes first. The recursion~\eqref{eqrec} comes from the fact that the root of one of the two sub-trees, parented by $v_0$, has to become sticky in order for the process to stop on the next step, and the paths of the process on these two sub-trees are independent of each other.

Note that  $\min\{a,b\}\le \eta(a,b) \le a+b-1$ and
\bn
\P(\eta(a,b)=\ell,\ \text{$a$ heads expired})&=\binom{\ell-1}{\ell-a}\frac{1}{2^\ell}, \quad a\le\ell\le a+b-1;
\\
\P(\eta(a,b)=\ell,\ \text{$b$ tails expired})&=\binom{\ell-1}{\ell-b}\frac{1}{2^\ell}, \quad b\le\ell\le a+b-1
\en
yielding
\bn
\P(\eta(a,b)=\ell)&=\left[
{\ell-1 \choose \ell-b}+
{\ell-1 \choose \ell-a}
\right]\frac{1}{2^\ell}, \quad \min\{a,b\}\le\ell\le a+b-1.
\en
with the convention ${x\choose y}=0$ if $y< 0$.

\vskip 5mm
\begin{figure}[!htb]\centering
\includegraphics[scale=0.5]{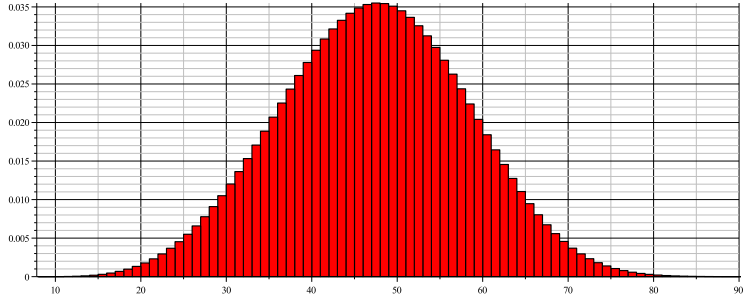}
\caption{Distribution of $\xi_8$}\label{fig:tree8}
\end{figure}
\vskip 5mm

Using~\eqref{eqrec} we can, in principle, get the distribution of~$\xi_K$ for any positive integer~$K$. For example, the distributions~of $\xi_4$ and $\xi_5$ are given in the following two tables:
\vskip 2mm
\begin{tabular}{|c|ccccc|}
\hline
$k$ & 4 & 5 & 6 & 7 & 8\\
\hline
$\P(\xi_4=k)$ & $\frac{1}{8}$& $\frac{1}{4}$& $\frac{5}{16}$& $\frac{15}{64}$& $\frac{5}{64}$
\\ \hline
\end{tabular}
\vskip 2mm
\noindent 

{\tiny
\vskip 2mm
\hskip -7mm
\begin{tabular}{|c|cccccccccccc|}
\hline
k&5&6&7&8&9&10&11&12&13&14&15&16 \\[1mm] \hline
$\P(\xi_4=k)$&$\frac{1}{64}$ &$\frac{3}{64}$ &$\frac{45}{512}$ &$\frac{535}{4096}$ &$\frac{1335}{8192}$ &$\frac{355}{2048}$ &$\frac{5115}{32768}$ &$\frac{30525}{262144}$ &$\frac{9075}{131072}$ &$\frac{32175}{1048576}$ &$\frac{75075}{8388608}$ &$\frac{10725}{8388608}$ 
\\[1mm] \hline
\end{tabular}
}
\vskip 2mm
For $\xi_8$ the distribution is shown on Figure~\ref{fig:tree8}. We have also computed
$$
\E \xi_3=3.5,\ 
\E \xi_4=5.89...,\ 
\E \xi_5\approx 9.82,\ 
\E \xi_6\approx 16.4,\ 
\E \xi_7\approx 27.6,\ 
\E \xi_8\approx 46.8.
$$
Our guess is that $\xi_K$,  appropriately scaled, is asymptotically normal for large $K$. Unfortunately, we do not have a proof of this fact, and leave this as a conjecture. One can also generalize the recursion~\eqref{eqrec} for regular $d$-ary trees with $d\ge 3$, but the formula quickly becomes quite messy and not so useful.

\subsection{Lower and upper bounds for $\xi_K$}
Here we deal with a regular $d-$ary tree again, dropping the restriction $d=2$.
Observe that the number of particles which get stuck on level $i$, $i=1,2,\dots,K-1$, (i.e.\ distance $i$ from the root), is at most $d^{i-1}$, since whenever a vertex $v$ becomes sticky, none of its sisters on the tree cannot be reached (if a new particle reaches the parent of $v$, it stops and becomes sticky). Therefore, since initially all the points on level $K$ are border points, a non-random upper bound on $\xi_K$ is given by
\begin{align}\label{eq:otherbound}
\xi_K\le \left[d^{K-2}+d^{K-3}+\dots+d+1\right]+1
\le d^{K-1} \left[\frac{1}{d-1}+d^{-(K-1)}\right]
\end{align}
where the last term ``$+1$" corresponds to the very last particle emitted at $v_0$ which immediately becomes sticky. The trivial lower  bound for $\xi_K$ is $K$, but we will show that with a high probability $\xi_K$ is in fact much larger.

Let $|v|=\dist(v,v_0)$ be the height of the particles on a tree, and let ${\sf Lev}_i=\{v:\ |v|=i\}$ be the set of $d^i$ vertices on level $i$. 
Let $\eta_i$ be the index of the particle which was first to get stuck on level $i$,
i.e.\  $\eta_i=\inf\{n\ge 1:\ v_n\in {\sf Lev}_{i}\}$. Trivially $\eta_{K-1}=1$; we will show that $\eta_{K-m}$ is quite large for $m\ge 2$. This will allow us to get the necessary bound as $$\eta_{K-1}<\eta_{K-2}<\dots<\eta_1 <\eta_0 \equiv \xi_K-1.$$

Fix some $m\ge 2$. Observe that for a vertex in ${\sf Lev}_{K-m}$ to become sticky, at least $m$ particles of random walk should pass through it on their way up. Since each vertex at level ${\sf Lev}_{K-m}$ is equally likely to be visited by the random walk (until there is at least one sticky particle at this level), the quantity $\eta_{K-m}$ is stochastically larger than $\zeta=\zeta_{K,m}$, the number of independent trials of a discrete uniform random variable with $A=d^{K-m}$ equally likely outcomes required to reach one of the $A$ outcomes {\em at least} $m$ times. Note that for $d=m=2$ and $A=365$. this is exactly the famous {\it birthday problem}; therefore
$$
\P(\zeta_{K,2}>t)=\prod_{i=1}^{t-1} \left(1-\frac{i}{d^{K-2}}\right)\sim \exp\left\{-\frac{t^2}{2\,d^{K-2}}\right\}
$$
and in particular if $h(t)$ is any positive function such that  $\frac{h(t)}{d^{t/2}}\downarrow 0$ then
$$
\P(\xi_K>h(K))\ge \P(\eta_{K-2}>h(K))\ge \P(\zeta_{K,2}>h(K))\to 1\text{ as }K\to \infty.
$$
For larger $m$, we do the following estimation. We have
\begin{align}\label{eqtreeA}
\P(\zeta_{K,m}\le t)&=
\P(\text{one of $A$ outcomes is reached $\ge m$ times during $t$ trials}) \nonumber
\\ \nonumber
&\le A\cdot \P(\text{outcome ``$1$'' reached at least $m$ times  during $t$ trials}) \\ \nonumber
&=A\cdot \sum_{i=m}^A {t\choose i} 
\frac{1}{A^i}\left(1-\frac{1}{A}\right)^{t-i}
\le \sum_{i=m}^A {t\choose i}  \frac{1}{A^{i-1}}
\sim {t\choose m}  \frac{1+o(1)}{A^{m-1}} 
\\
&=\frac{t^m+o(t^m)}{m! A^{m-1}} 
\end{align}
if $t\ll A$ and 
$m\ll t$.
By Stirling's formula, the logarithm of the RHS of~\eqref{eqtreeA} is approximately
\begin{align}\label{eqtreeAA}
&m\ln t - (m\ln m - m+\OO(\ln m)) -(K-m)(m-1)\ln d \nonumber
\\
&=m\ln t  -K(m-1)\ln d+m^2 \ln d +\OO(m\ln m)
\end{align}
We want this quantity to be negative and to go to $-\infty$, but preferably slowly. Equating the RHS of~\eqref{eqtreeAA} (but the $\OO(\cdot)$ term) to $0$ gives
\begin{align*}
t=d^{K(1-1/m)-m}=d^{K -\left[m+\frac{K}{m}\right]}
\end{align*}
and substituting this into the LHS of~\eqref{eqtreeAA}
we get
$$
- m(\ln m+\ln d  -1)+\OO(\ln m) 
$$
Since we want $t$ to be as large as possible, we  choose an integer $m=\sqrt{K}+z$ such that $|z|\le 1/2$.
Then, indeed,  $t=d^{K-K/m-m}\ll A=d^{K-m}$ and $m\ll t$; moreover, the RHS of~\eqref{eqtreeAA} becomes $-\left(\frac 12+o(1)\right)\sqrt K \ln (K) \to-\infty$ as $K\to\infty$.
Hence, taking into account that $m+\frac{K}{m}\le 2\sqrt{K} +\OO(1/\sqrt{K})$, we get
\begin{align}\label{eqtreeAAA}
\P\left(\zeta_{K,m}>d^{K-2\sqrt K-\OO(K^{-1/2})}\right)\to 1 \quad\text{as }K\to\infty.
\end{align}
Since $\eta_K>\zeta_{K,m}$, combining with~\eqref{eq:otherbound}, we have the following statement.
\begin{thm}
With probability at least $1-K^{-\left(\frac 12+o(1)\right) \sqrt{K}}$ we have for $d\ge 2$
$$
1+\frac{\ln (d-1)}{\ln d}-o(1)\le K-\log_d \xi_K \le 2\sqrt{K}+\OO\left(K^{-1/2}\right).
$$
\end{thm}

\section{Two-dimension aggregation on $\Z^2$}\label{Sec:2-dim-lattice}

\begin{figure}[!htb]\centering
\includegraphics[scale=.3]{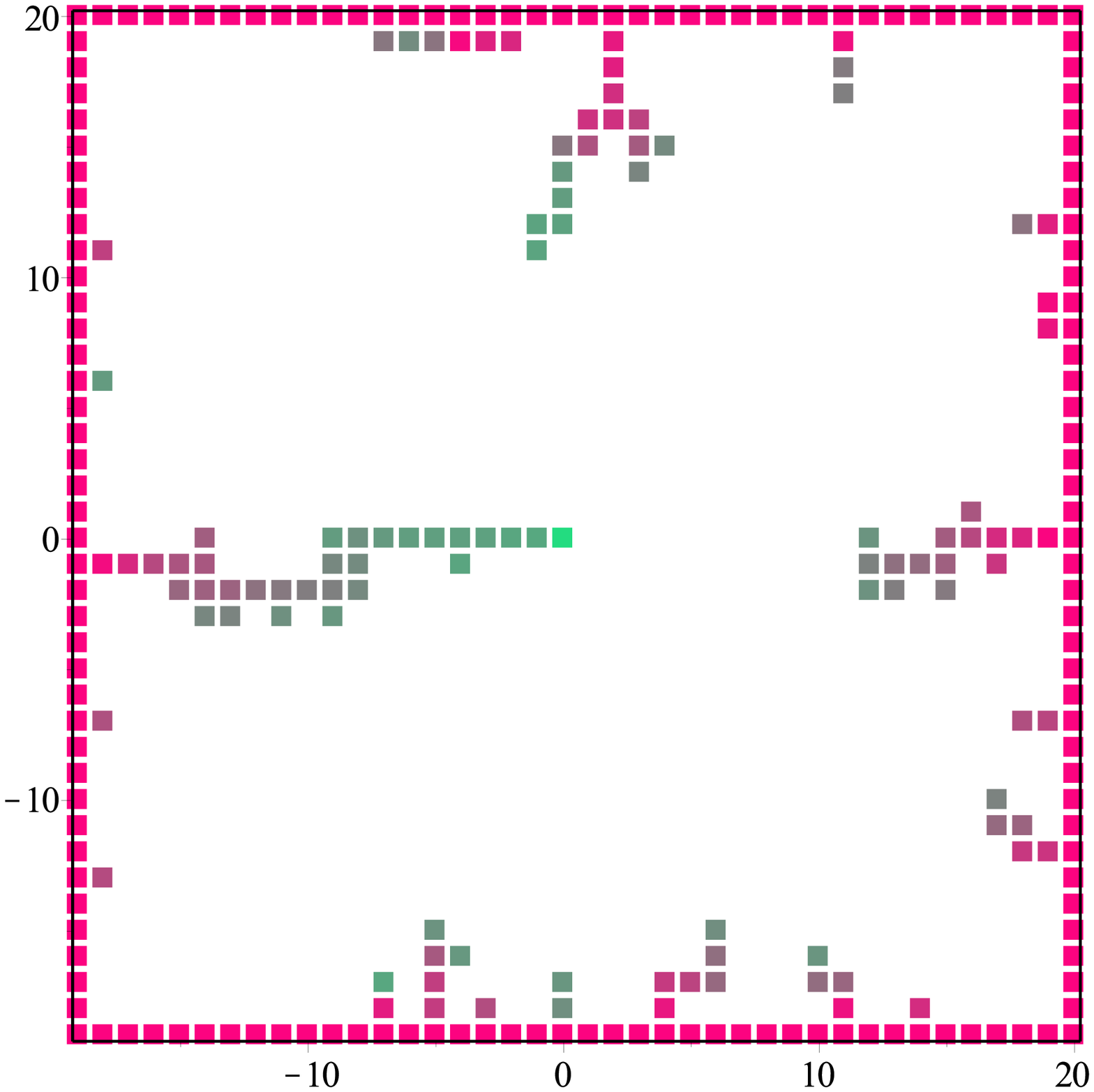}
\hskip 1cm
\includegraphics[scale=.3]{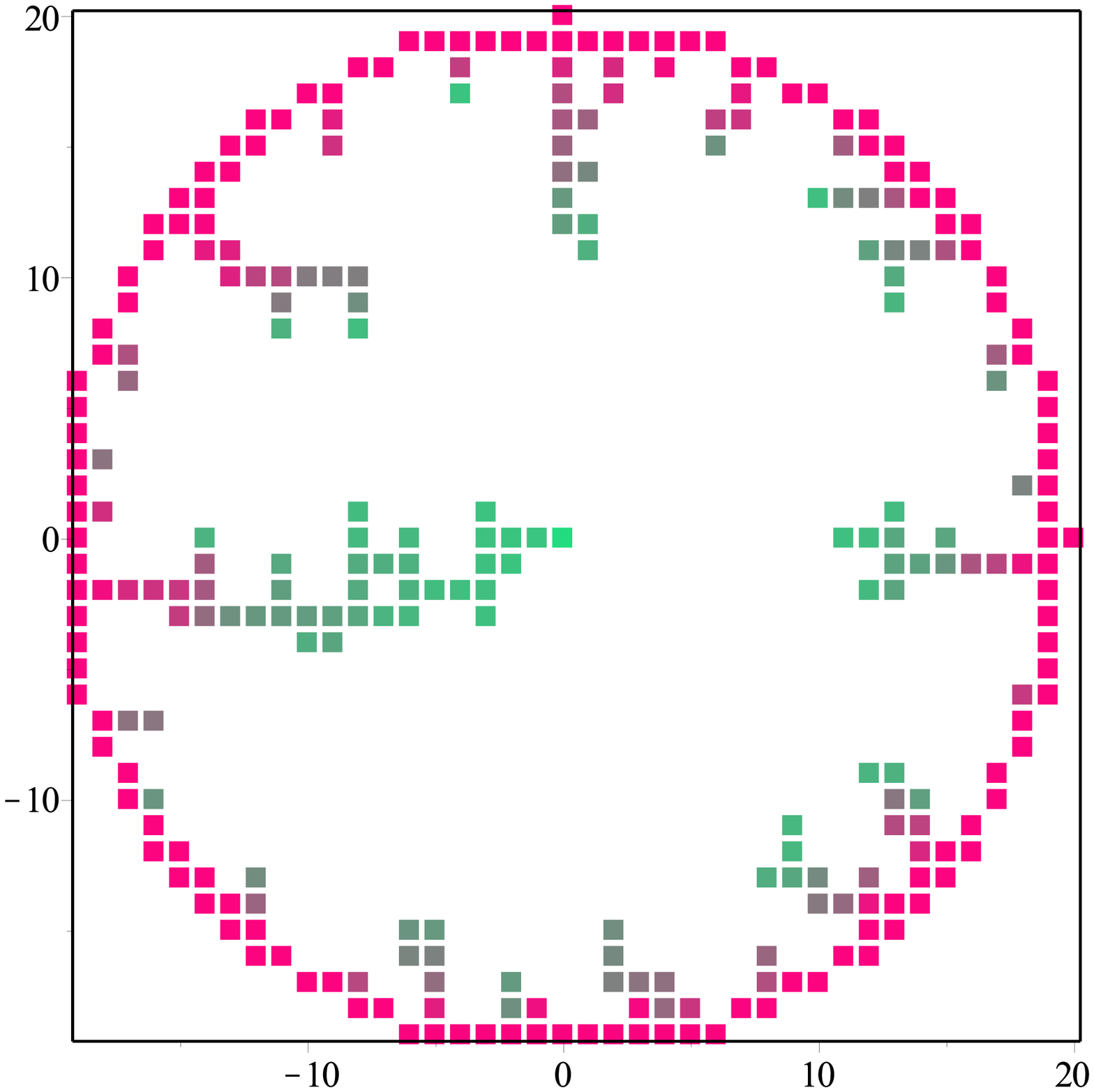}
\caption{Aggregation process on $\Z^2$ with square (a) and circle (b) borders}\label{fig:digraph}
\end{figure}

Let the graph $G$ be a box $[-N,\dots,N]^2\subset \Z^2$ with the origin $v_0=(0,0)$. We can define the model in two alternatives ways:
\begin{itemize}
\item[(a)]
(box model) ${\sf B}=\{(x,y)\in G:\ |x|=N\text{ or } |y|=N\}$ is the border of the box $[-N,\dots,N]^2$;
\item[(b)]
(disc model) ${\sf B}=\{(x,y)\in G:\ \sqrt{x^2+y^2}\ge N-1\}$, i.e.\ $G$ can be viewed as the disc of radius $N$ and the ``sticky border" is the circumference. 
\end{itemize} 
Figure~\ref{fig:digraph} shows the aggregation process at the time when the process has stopped, compare with illustrations in~\cite{LE&Pe07}. In what follows, we study only case (b).

Let $\xi_N$, as before, denote the number of emitted particles until the origin becomes part of the border. It is trivial that $N\le \xi_N\le N^2$, however we want to get a finer asymptotic of $\xi_N$; we conjecture that $\xi_N\sim N^{\a}$ where $\alpha\approx 1.7$ (see also~\cite{LE&Pe07}), however we believe it is a very hard problem. We have managed, though, to show that $\xi_N$ is at least of higher order than $N^{4/3-\eps}$, please see Theorem~\ref{thm:N43} below.

\subsection{Lower bound for the BA model on a disc }\label{subsecbox2}
\begin{thm}\label{thm:N43}
For every $\eps>0$ we have
$$
\P\left(\xi_N<  N^{\frac{4}{3}-\eps} \right)\to 0
$$ 
as $N\to\infty$.
\end{thm}
The proof will be based on obtaining detailed bounds for the DLA of~\cite{Ke87,Ke87_hit} via strengthening of the result of the main theorem in~\cite{Ke87}. In accordance with notations of these papers, let $B\subset \Z^2$ be a finite connected subset, $\partial B$ be the set of points adjacent to $B$, $S_n$ is a simple symmetric random walk on $\Z^2$ with $S_0=x$, $\tau=\inf\{n:\ S_n\in B\}$ the a.s.\ finite hitting time of $B$, and $S_{\tau}$ is the point where the walk hits $B$ for the first time. Let also for $x\notin B$
$$
H(x,y)=\P_x(S_\tau=y),\quad
\mu(y)=\lim_{|x|\to\infty} H(x,y).
$$
The latter limit exists and satisfies $\sum_{y\in B}\mu(y)=1$ according to~\cite{Spit}, Theorem~14.1.
Suppose that $B$ contains the origin. Let $r=r_B=\max_{x\in B}|x|$ be the ``radius'' of $B$.
Kesten~\cite{Ke87_hit} showed that
$$
\mu(y)\le \frac{\sf const}{\sqrt{r_B}}
$$
where the constant does not depend on $B$.
We want first to generalize this result for {\em finite} starting point $x$.
\begin{prop}\label{propKK}
There is a constant $C>0$ not depending on anything, such that
$$
H(x,y)\le \frac{C}{\sqrt r_B}
$$
for all $B$ containing the origin and $x$ satisfying  $|x|\ge r_B^{3/2}\ln r_B$.
\end{prop}
\begin{proof}
Throughout the proof we fix the set $B$ and will write $r=r_B$ for simplicity. In accordance with the notations in~\cite{Spit}, let  $P_n(x,y)$ denote the $n$-step transition probability from $x$ to $y$ for a SRW on $\Z^2$, $$
G_n(x,y)=\sum_{k=0}^n P_k(x,y),\quad 
A(x,y)=\lim_{n\to\infty}\left[G_n(0,0)-G_n(x,y)\right],
$$
that is, $G_n(x,y)$ is the $n-$step Green function (see~\cite{Spit} Defintions~1.4 and~11.1).
We use the representation for $H(x,y)$ from~\cite{Spit}, formula~(14.1), and the proof of Theorem~14.1 there which states that 
\begin{align*}
H(x,y)-\mu(y)&=\sum_{t\in B} A(x,t) \left[\Pi(t,y)-1_{t=y}\right]
\\
&=\sum_{t\in B} \left[A(x,t)-A(x,0)\right] \Pi(t,y)
-\left[A(x,y)-A(x,0)\right]
\end{align*}
where $\Pi(t,y)\ge 0$ denotes the probability that the first hit into $B$ starting from $t$ will be at point $y\in B$ and satisfies $\sum_{t\in B} \Pi(t,y)=1$ by Lemma~11.2 in~\cite{Spit}.
 Hence
\begin{align}\label{eqHmu}
|H(x,y)-\mu(y)|& \le 
\left| A(x,y)-A(x,0)\right|+\sum_{t\in B} \left| A(x,t)-A(x,0)\right| \Pi(t,y)
\nonumber \\ &\le \nonumber
\max_{t\in B} \left|A(x,t)-A(x,0)\right| 
\left(1+\sum_{t\in B} \Pi(t,y)\right)
\nonumber \\ &
=2\max_{t\in B} \left|A(x,t)-A(x,0)\right| .
\end{align}

Next we need to estimate how quickly the difference between~$A(x,t)$ and~$A(x,0)$ converges to $0$.
From the proof of Proposition~12.2 in~\cite{Spit} and the translation invariance of the walk it follows that
\begin{align*}
A(x,t)-A(x,0)&=\lim_{n\to\infty} G_n(x,0)-G_n(x,t)
=\lim_{n\to\infty} G_n(x,0)-G_n(x-t,0)
\\
&=a(x-t)-a(x)=\frac{1}{(2\pi)^2}
\int_{-\pi}^{\pi}\int_{-\pi}^{\pi} 
Q\cdot R \, \rmd\theta_1 \rmd\theta_2 
=:(*)
\end{align*}
where
\begin{align*}
Q&=e^{i (x_1\theta_1+x_2\theta_2)}, 
\quad
R=\frac{1-e^{-i(t_1\theta_1+t_2\theta_2)}}{1-\phi}, 
\ \text{ and }\ 
\phi=\frac{\cos(\theta_1)+\cos(\theta_2)}{2}
\end{align*}
is the characteristic function of the walk. While it follows from~\cite{Spit} that the integral $(*)\to 0$ as $|x|\to\infty$, we still have to estimate the speed of this convergence.
Assume w.l.o.g.\ that $|x_2|\ge |x_1|$ so that $|x_2|\ge |x|/2$. Split the area of integration $[-\pi,\pi]^2$ into two parts: where $|\theta_1|<\eps$ and the remaining part, and write $(*)=(I)+(II)$ where (I) is the integral over the first area and (II) is the integral of the remaining area.

First, we obtain two useful inequalities:
\begin{align*}
\left|1-e^{-i(t_1\theta_1+t_2\theta_2)}\right|&=
2\left|\sin\left(\frac{t_1\theta_1+t_2\theta_2}{2}\right)\right|
\le |t_1\theta_1+t_2\theta_2|
\\
&\le \sqrt{2}\max(|t_1|,|t_2|)\sqrt{\theta_1^2+\theta_2^2}
\le 2r \,\sqrt{\theta_1^2+\theta_2^2}
\end{align*}
and if $|u|\le \pi$ then 
\begin{align*}
1-\cos(u)&\ge \frac{u^2}{5}
\ \Longrightarrow 2-\cos(\theta_1)-\cos(\theta_2)
\ge \frac15\left(\theta_1^2+\theta_2^2\right).
\end{align*}
Integrating $Q\cdot R$ by parts w.r.t.~$\theta_2$
for $|\theta_1|>\eps$ we get
\begin{align*} 
(I)&=\frac{1}{(2\pi)^2}
\int_{\eps<|\theta_1|<\pi} \rmd \theta_1
\int_{-\pi}^\pi \rmd \theta_2  \ Q \cdot R\\
&=
\frac1{(2\pi)^2}\int_{\eps<|\theta_1|<\pi}\rmd \theta_1
\left[ 
\frac{R\cdot e^{i(x_1\theta_1+x_2\theta_2)}}{ix_2}\,
 \Big|_{\theta_2=-\pi}^{\pi} 
-\int_{-\pi}^{\pi} 
\frac{e^{i(x_1\theta_1+x_2\theta_2)}}{ix_2}\, \frac{\partial R}{\partial \theta_2}
\rmd \theta_2\right]
\\
&=
\frac{i}{(2\pi)^2 \, x_2}\int_{\eps<|\theta_1|<\pi}\rmd \theta_1
\int_{-\pi}^{\pi} 
Q\, \frac{\partial R}{\partial \theta_2} \rmd \theta_2
\end{align*}
using the fact that $e^{i\pi x_2}=e^{-i\pi x_2}=0$ as $x_2\in\Z$. Since
\begin{align*}
\frac{\partial R}{\partial \theta_2}
&=\frac{
i t_2 e^{-i(\theta_1 t_1+\theta_2 t_2)}
\left[2-\cos\theta_1 -\cos\theta_2\right]
+
\left[e^{-i(\theta_1 t_1+\theta_2 t_2)}-1\right]
\sin \theta_2
}{(2-\cos\theta_1 -\cos\theta_2)^2}
\end{align*}
and $|\sin\theta_2|\le |\theta_2|$ we conclude that
\begin{align*}
(2\pi)^2 x_2 \, |(I)|&\le
\int_{\eps<|\theta_1|<\pi}\rmd \theta_1
\int_{-\pi}^{\pi} 
\frac{4|t_2|}{2-\cos\theta_1 -\cos\theta_2}
\rmd \theta_2
\\
&+
\int_{\eps<|\theta_1|<\pi}\rmd \theta_1
\int_{-\pi}^{\pi} 
\frac{  
2r|\theta_2|\,\sqrt{\theta_1^2+\theta_2^2}
}{(2-\cos\theta_1 -\cos\theta_2)^2}\rmd \theta_2
=(Ia)+(Ib).
\end{align*}
Let us write $r=r_B$ to simplify the notations. Since $|t_2|\le r$ we have
\begin{align*}
(Ia)&\le 20r
\int_{\eps<|\theta_1|<\pi}\rmd \theta_1
\int_{-\pi}^{\pi} 
\frac{\rmd \theta_2}{\theta_1^2+\theta_2^2}
< 20r \iint_{\eps^2<\theta_1^2+\theta_2^2<4\pi^2} \frac{\rmd \theta_1\rmd \theta_2}{\theta_1^2+\theta_2^2}
\\
&=20r\int_{\eps}^{2\pi} \frac{\rmd \rho}{\rho}
\int_{-\pi}^{\pi} \rmd \varphi 
=40\pi r \ln(2\pi \eps^{-1})
\end{align*}
by switching to polar coordinates. Similarly,
\begin{align*}
(Ib)&\le 50r
\int_{\eps<|\theta_1|<\pi}\rmd \theta_1
\int_{-\pi}^{\pi} 
\frac{|\theta_2|\rmd\theta_2}{(\theta_1^2+\theta_2^2)^{3/2}}
=
200r
\int_{\eps}^\pi
\frac{1}{\theta_1}-\frac{1}{\sqrt{\pi^2+\theta_1^2}}
\rmd \theta_1
\\ &=200 r (1+o(1))\, \ln(\eps^{-1}).
\end{align*}
Consequently,
$$
|(I)|\le \left|\frac{(Ia)+(Ib)}{4\pi^2 x_2}\right|\le  \frac{50+10\pi+o(1)}{\pi^2}\cdot \frac{2r}{|x|}\cdot \ln (\eps^{-1}).
$$
\vskip 5mm
On the other hand, for $|\theta_1|\le\eps$ we have
\begin{align*}
|(II)|&\le \frac{2}{(2\pi)^2}\int\int_{|\theta_1|\le \eps} \frac{|1-e^{-i(t_1\theta_1+t_2\theta_2)}|}{2-\cos(\theta_1)-\cos(\theta_2)}\rmd \theta 
\le 
 \frac{5r}{\pi^2}
 \int\int_{|\theta_1|\le \eps} \frac{\sqrt{\theta_1^2+\theta_2^2}}{\theta_1^2+\theta_2^2}\rmd \theta \\
&=
\frac{5r}{\pi^2} 
\left[
4 \pi  \ln\left( 
\frac {\epsilon}{\pi}
+\sqrt{1+\frac {\epsilon^2}{\pi^2}} 
 \right) +2\epsilon \ln  \left( 
 \frac{\sqrt {\pi^2+{\epsilon}^2}+\pi}
 { \sqrt {{\pi }^{2}+{\epsilon}^{2}}-\pi}
 \right) 
\right]
\\
&=\frac{20r}{\pi^2} \left[\ln(2\pi e)+\ln \left(\eps^{-1}\right)\right]\cdot \epsilon	 +\OO \left( \eps^3 \right) =\OO \left( r\eps \ln \left(\eps^{-1}\right) \right) 
\end{align*}

Therefore, by setting $\eps=\frac{1}{r^{3/2} (\ln r)^2}$ we get
\begin{align*}
|(A(x,t)-A(x,0)|&=|(*)|
\le |(I)|+|(II)|
\le {\rm const} \left[\frac{r \ln(\eps^{-1})}{|x|}+r \eps \ln(\eps^{-1})\right]
\\ &=
{\rm const} \left[\frac{r\ln\left(r^{3/2} (\ln r)^2\right)}{|x|}
+\frac{
\ln\left(r^{3/2} (\ln r)^2\right)
}{(\ln r)^2\ \sqrt r } 
\right]
\\
 &=
\frac{\rm const}{2} \left[r\frac{3 \ln r + 4\ln\ln r}{|x|}
+\frac{
3 \ln r + 4\ln\ln r
}{(\ln r)^2 \ \sqrt r } 
\right]
\\ & \le \frac{3 \cdot {\sf const}}{2} \cdot\frac{1}{\sqrt{r}}+ o\left(\frac{1}{\sqrt{r}}\right).
\end{align*}
Now the result follows from~\eqref{eqHmu}, Theorem from~\cite{Ke87_hit}, and the condition of our theorem about~$|x|$.
\end{proof}

Next we want to strengthen Kesten's result from~\cite{Ke87}, where he studied the following model.
Suppose that the initial sticky particle is located at the origin ${\bf 0}\in\Z^2$, and the particles emitted at infinity (for more rigorous definition please see~\cite{Ke87}). Let $r(n)$ be the diameter of the aggregate $A_n$ of $n$ particles. Kesten in~\cite{Ke87} showed that a.s.\ $r(n)\le C_2 n^{2/3}$ for a fixed constant $C_2>0$ and all but finitely many $n$. We will get a more precise estimate for all $n$ even in the case where the particles are emitted {\em not at the infinity} but at some point, sufficiently remote from the origin.

\begin{prop}[Strengthening Kesten's theorem for $\Z^2$]\label{propsuperKesten}
Consider the above  model with the exception that the particles are emitted from a fixed finite point $z\in\Z^2$. 
Then there are constants $C_4,C_5,n_0>0$ not depending on anything, such that for all  $n\ge n_0$  satisfying $n^{3/2}\ln n\le |z|$ 
we have
$$
\P(r(n)> C_4 n^{2/3} ) \le e^{-C_5\sqrt n}.
$$
\end{prop}
\begin{proof}
We assume that $n=n_i=2^{2i}$ for some positive integer $i$; if this is not the case then we can always find such $i$ that $n_{i-1}\le n< n_{i}$ and since $n_i/4\le n\le n_i$ and $A_n\subseteq A_{n_i}$ the result will follow.

For $\sqrt{n}=2^i$ we have a trivial bound $r(2^i)\le 2^i$. Now for $k=i,i+1,...,2i$ we repeat Kesten's argument. Note here that our Proposition~\ref{propKK} together with the trivial bound $r(A_n)\le n$ imply that the inequality~(8) from~\cite{Ke87} still holds, possibly with a different constant; that is, the probability that the particle get adsorbed at a specific point of $A_n$ is bounded by 
$\displaystyle \frac{C_6/16}{\sqrt{r(n)}}$ for some $C_6>0$.
Then the collection of inequalities~(9) in~\cite{Ke87}, that is
\begin{align}\label{eqKe8}
r(2^{k+1})-r(l)\le \frac{C_6 2^k}{\sqrt{r(l)}}+2^{k/2}
\quad \text{for all } 2^{k}\le l\le 2^{k+1}
\end{align}
holds with probability at least $1-\nu(k)$ where
$$
\nu(k)=4\pi^2 2^{k+1} \left(\frac{e}{4}\right)^{2^{k/2}}
\times 2^k \le \g^{2^{k/2}}
$$
for some $\g<1$ and all $k$ larger than some non-random $k_0$ (see equation~(18) in~\cite{Ke87}). 

From now on assume that $i>k_0$, that is, $n\ge 4^{k_0}$. Then~\eqref{eqKe8} holds for all $k=i,i+1,\dots,2i-1$ with probability exceeding
\begin{align}\label{eqgam}
1-\sum_{k=i}^{2i} \g^{2^{k/2}}
> 1-\sum_{m=2^{i/2}}^{\infty} \g^{m}
= 1- \frac{\g^{2^{i/2}}}{1-\g}=
1- \frac{\g^{\sqrt{n}}}{1-\g}
\end{align}

Next, suppose that inequalities~\eqref{eqKe8} indeed hold for $2^k\le l\le 2^{k+1}$. If for $l=2^k$ we have $r(l)\ge (C_6 2^{k})^{2/3}$ then 
$$
r(2^{k+1})-r(2^k)\le  (C_6 2^{k})^{2/3} +2^{k/2}.
$$
If the above inequality {\it does not} hold, then either for all $l\in [2^k,2^{k+1}]$ we have $r(l)<(C_6 2^{k})^{2/3}$ and thus $r(2^{k+1})<(C_6 2^{k})^{2/3}$ as well, or there is some $l_*\in [2^k,2^{k+1}]$ such that $\left| r(l_*)-(C_6 2^{k})^{2/3}\right|\le 1$. In the latter case
\begin{align*}
r(2^{k+1})&\le r(l_*)+ \frac{C_6 2^k}{\sqrt{r(l_*)}} +2^{k/2}
\\ &
\le \left[(C_6 2^{k})^{2/3}+1\right]+
\frac{C_6 2^k}{\sqrt{(C_6 2^{k})^{2/3}-1}}+2^{k/2}
\le 3(C_6 2^{k})^{2/3}+2^{k/2}.
\end{align*}
Combining the inequalities involving $r(2^{k+1})$ we conclude that
$$
r(2^{k+1})\le  3\,C_6^{2/3}\, 2^{\frac{2k}3} +2^{\frac k2}+r(2^k)
\le 4\, C_6^{2/3}\,  2^{\frac{2k}3} +r(2^k)
$$
if $k$ is not too small. Summing this up for $k=i,i+1,\dots,2i-1$ we finally get
$$
r(n)=r(2^{2i})\le r(2^i)+ 4\,C_6^{2/3}\, \sum_{k=i}^{2i-1}2^{\frac{2k}3}
\le
2^i+\frac{4\,C_6^{2/3}}{1-2^{2/3}}\,2^{\frac{4i}3}
\le C_4\, 2^{\frac{4i}{3}} =C_4\, \cdot n^{\frac{2}{3}}
$$
with probability exceeding the quantity in the RHS of~\eqref{eqgam}.
\end{proof}

\begin{corollary}\label{corsuperK}
Let $\a>0$ be small. 
Consider again the  model from Proposition~\ref{propsuperKesten}, with the same $z$.
Then there are constants $C_4,C_5,n_0>0$ depending only on $\a$, such that for all  $n\ge n_0$  satisfying $n\le |z|^{1-\a}$
we have
$$
\P(r(n)> C_4 n^{2/3} ) \le e^{-C_5\sqrt n}.
$$
\end{corollary}
\begin{proof}
The crucial point 
in the proof of Proposition~\ref{propsuperKesten} where we used the fact that $|z|\ge n^{3/2}\ln n$ was that we can apply Proposition~\ref{propKK} only as long as the set $B$ containing the origin has the radius $r=r(B)$ satisfying $r^{3/2} \ln r\le |z|$. The estimate $r(A_n)\le n$ which we used in the beginning of the proof of Proposition~\ref{propsuperKesten} is, however, too crude, as we know that $A_n$ does not grow that fast with very high probability. Therefore, one can repeat the arguments of this proposition almost verbatim, estimating the probabilities {\em conditioned} on the past behaviour of the adsorption process not to grow faster then $s^{2/3}$ by time $s$ for $s\le n$, so that in particular $r(n)\le C_4 n^{2/3}$ and thus 
$r(n)^{3/2} \ln r(n) \le \OO(n \ln n)\ll |z|$ as it would be required by the conditions of Proposition~\ref{propsuperKesten}.
\end{proof}

Now we present  the proof of the main result, based on estimation of crossing times of the sequence of rings separating the border from the origin.

\begin{proof}[Proof of Theorem~\ref{thm:N43}.]
Let $B(R)=\{v\in Z^2:\ |v|\le R\}$ be the set of points in $\Z^2$ inside the circle of radius $R$; fix some very small positive $\d$ such that 
\begin{align}
\label{epsdel}
\frac{\d}{18-6\d}<\frac{\eps}2
\end{align}
and let $r_k=k^{3-\d}$, $k\in \Z_+$, and for now assume that $N=r_z$ for some positive integer~$z$. Consider the rings $R_k=B(r_k)\setminus B(r_{k-1})$.
Observe that the width of $R_k$ is $\cong (3-\d)k^{2-\d}$ and thus it is larger than $w_k:=k^{2-\d}$.

Next step is to show that with high probability the number of particles $\zeta_k$ required to cross $R_k$ and come to the next ring $R_{k-1}$ (even if there is more than one ``arm'', that is, a connected component of stuck particles) is at least of order $ w_k^{3/2}= k^{3-\frac{3}{2}\d}$. with high probability.

Consider our adsorption process from the moment when some particle gets adsorbed in $R_k$ for the first time. Let $\partial B(r_k)$ be the set of vertices where this could have happened, namely
$$
\partial B(r_k)=\{v\in B(r_k): \ \exists u\notin B(r_k),\ u\sim v\}
$$
``the internal border'' of $B(r_k)$. Note that ${\sf card}(\partial B(r_k))\le 8 r_k$. 

Let us arbitrarily index points of $\partial B(r_k)$ as $v_j$, $j=1,2,\dots,{\sf card}(\partial B(r_k))$. For each $v_j\in \partial B(r_k)$ construct the corresponding ``DLA arm'' $A^j\subseteq R_k$ as follows.
Initially all $A^j$ are empty sets. Whenever a particle gets adsorbed in a point $v_j\in \partial B(r_k)$  set $A^j=\{v_j\}$. If a particle gets adsorbed at some previously empty point $v\in R_k\setminus \partial B(r_k)$, then for every $u$ such that $u\sim v$ and every $j$ such that $u\in A^j$, attach $v$ to $A^j$, i.e., $A^j\to A^j \bigcup \{v\}$.
\footnote{Observe that as a result point $v\in R_k$ can {\it simultaneously} join a number of different ``arms''.}
 Finally, if the particle gets adsorbed outside of $R_k$, do not change any of the arms.

Formally, let $t$ be the index of the particle emitted from the origin counting from the first time a particle got adsorbed in $R_k$ at some point $v\in \partial B(r_k)$. Then
$$
A^j_1=\emptyset \text{ for all but one $i$ for which }
A^i_1=\{v\}.
$$
Now recursively define $A^j_t$, $t=2,3,\dots$, as follows: for each $j$
\begin{align*}
A^j_{t+1}=\begin{cases}
A^j_t \cup \{v\}
&\text{if the $t+1$-st particle got adsorbed at $v$ such that } v\sim A^{j}_t\\
A^j_t&\text{otherwise}.
\end{cases}
\end{align*}
It is clear from the construction  that for any arm $A^j$, the probability to get adsorbed near any of its points is smaller than the corresponding probability for the process described in Proposition~\ref{propsuperKesten} in particular, the number  of particles in $A^j$ might grow slower than the number of emitted from the origin particles, i.e.\ ${\sf card} \left(A^j_t\right)\le t$ (unlike the Kesten's DLA model where ${\sf card}(A(t))=t$\,).

Set $n=n_k=(w_k/C_4)^{3/2}=k^{3-1.5\d}\, C_4^{-1.5}$. Then 
$$
n =C_4^{-1.5} \left[k^{3-\d}\right]^{1-\frac{\d}{6-2\d}}
=\OO\left( r_{k-1}^{1-\frac{\d}{6-2\d}}\right)
$$
so we can apply Corollary~\ref{corsuperK} with $\a\in\left(0,\frac{\d}{6-2\d}\right)$ to show that  after $n$ particles were adsorbed inside of $R_k$, for each arm $A^j$ we have
$$
\P(r\left(A^j_n\right)>w_k) = \P(r\left(A^j_n\right)>C_4 n^{2/3}) \le e^{-C_5\sqrt{n}}
$$
where $r(A^j_n)$ denotes the diameter of the arm $A^j_n$. For a path of sticky particles to cross the ring $R_k$ it is necessary that the diameter of at least one of the arms exceeds $w_k$. The probability that it took no more than $n_k$ particle to cross $R_k$ is
\begin{align*}
\P(\zeta_k\le n_k)&\le \P\left(\bigcup_{j}\left\{r(A^j)>w_k\right\}\right)
\le \sum_{j} \P(r(A^j)>w_k)
\\
&
\le 8 r_k e^{-C_5\sqrt{n_k}}
= 8k^{3-\d} e^{-C_5\sqrt{k^{3-1.5\d}\, C_4^{-1.5}}}
\le e^{-k}
\end{align*}
(since there are at most $8r_k$ such arms) for all $k$ large enough.

Fix an arbitrary $\eps'>0$ and choose $k_0$ so large that $\sum_{k=k_0}^{\infty} e^{-k}<\eps'$.
As a result, with probability 
$$
1-\sum_{k=k_0}^z e^{-k}>1-\eps'
$$ 
the number of particles required to form a path that crosses all the rings $R_z$, $R_{z-1}$, $\dots$, $R_{k_0+1}$, $R_{k_0}$ is no less than
$$
\sum_{k=k_0}^{z}  n_k
=\frac{1}{C_4^{1.5}} 
\sum_{k=k_0}^{z} k^{3-1.5\d} 
\ge 
\frac{z^{4-1.5\d}}{5 C_4^{1.5} } =\frac{ N^{\frac{4}{3}-\frac{\d}{18-6\d}}}
{5 C_4^{1.5} }
\ge \frac{ N^{\frac{4}{3}-\frac{\eps}2}}
{5 C_4^{1.5} }
\ge N^{\frac{4}{3}-\eps}
$$
(see~\eqref{epsdel}) provided $z\ge 2 k_0$ and $N$ is sufficiently large. This implies that
$$
\P\left(\xi_N\le  N^{\frac{4}{3}-\eps}\right) \le \eps'.
$$
Finally, if $N\ne z^{3-\d}$ for an integer $z$, we can always find $N'$ such that $N'=z^{3-\d}$ and  $N/2<N'\le N$ and apply the argument for the rings starting with $N'$.
\end{proof}




\subsection{Comb lattice}
The comb lattice $G$ is the graph whose vertices coincides with the vertices of $\Z^2$, however, all the horizontal edges are removed except those lying on the horizontal axes. Thus, a simple random walk located at point $(x,y)\in\Z^2$ goes only up or down ($y\pm 1$) with probability $\frac 12$, unless $y=0$ in which case either of the coordinates can decrease or increase, all with probability $\frac 14$.

Suppose the origin is $v_0=(0,0)$ and the initial sticky border consists of two horizontal lines located at distance $N$ from the horizontal axes, i.e.\  
${\sf B}=\{(x,y)\in G:\  |y|=N\}$. As before, let $\xi=\xi_N$ denote the number of particles to be emitted from the origin before the origin becomes sticky. 

\setlength{\unitlength}{1cm}
\begin{picture}(16,9)(2,0)
\put(1.5,     4){\vector(1,0){13}}
\put(2.5,     4){\vector(-1,0){1}}
\multiput(2,7)(1, 0){13}{\circle*{0.15}}
\multiput(2,6)(1, 0){13}{\circle*{0.15}}
\multiput(2,5)(1, 0){13}{\circle*{0.15}}
\multiput(2,4)(1, 0){13}{\circle*{0.15}}
\multiput(2,3)(1, 0){13}{\circle*{0.15}}
\multiput(2,2)(1, 0){13}{\circle*{0.15}}
\multiput(2,1)(1, 0){13}{\circle*{0.15}}
\multiput(2,1)(1,0){13}{\line(0,1){6}}
\put(8,7.1){\oval(13,0.6)}
\put(8,1.1){\oval(13,0.6)}
\put(6,0){Comb lattice, $N=3$}
\put(8.1,   4.3){$v_0$}
\put(15,1){${\sf B}$}
\put(15,7){${\sf B}$}
\end{picture}

\begin{thm}\label{thm:comb}
For some $0<c_1<c_2$ 
$$
\P(c_1 N^{3/2}\le \xi_N \le c_2 N^{3/2})\to 1
$$
as $N\to\infty$.
\end{thm}
The proof of Theorem~\ref{thm:comb} will immediately follow from Lemmas~\ref{lem:comblower} and~\ref{lem:combupper} below.

Consider the column $j$, $j\in\Z$, that is, the set of points $K_j=K_j^+\cup K_j^-$ where
$$
K_j^+=\{(j,y),\ y=0,1,\dots,N\}, \quad
K_j^-=\{(j,y),\ y=0,-1,\dots,-N\}.
$$
Both columns are eventually being filled with sticky particles; let $h_j^+(m)$ be the distance from $(j,0)$ to the closest sticky particle in $K_j^+$ at the time when the $m$-th particle is being emitted from the origin; similarly define $h_j^-(m)$  for $K_j^-$.

Suppose that when the $m$-th particle is emitted, all $h_j^{\pm}(m)\ge N/2$ for all $j\in\Z$. Consider the embedded random walk restricted to the horizontal axes ($y=0$), and denote it by $W_n\in\Z$, $W_0=0$. Eventually the particle gets stuck during an excursion to one of the columns when it reaches the sticky border there; thus this walk is defined only until some random stopping time $\tau=\tau(m)$, and the column in which it gets stuck is either $K_{x(m)}^+$ or $K_{x(m)}^-$ where $x(m)=W_{\tau(m)}$. We shall say  then that the walk $W_n$ ``dies" at time $\tau(m)$.

It is easy to see that up to time~$\tau$ the process $W_n$ is essentially a simple random walk on~$\Z^1$. From elementary calculations, given that~$W_n=j$, the probability  the walk dies before ever visiting $j\pm 1$  again (which means it reaches a sticky border in either $K_j^\pm$ before departing for $j-1$ or $j+1$) is given by
\begin{align}\label{eq:dienow}
q_{j,n}=q_{j,n}(m):=\frac{1/2}{1/2+\frac{1}{1/h^+ +1/h^-}}
=\frac{h^+ + h^-}{2h^+  h^- + h^+ + h^-}
\end{align}
where we omitted the subscript $j$ for simplicity (one can use e.g.\ electrical networks method, see~\cite{DS}.)
Consequently, if all $h_j^\pm \ge N/2$ (and by the initial conditions we know that $h_j^\pm\le N-1$) then 
\begin{align}\label{eq:qbounds}
\frac{1}{N}\le q_{j,n} \le \frac{2}{N+2} <\frac 2N.
\end{align}
Under the above assumption $\min_j h_j^\pm(j)\ge N/2$ we can compute the probability that the $m$-th particle eventually gets stuck at point~$j$ by
$$
p_j(m):=\sum_{s=0}^{\infty} \sum_{l\in{\cal L}_{s,j}} \frac{1}{2^s}
\cdot (1-q_{l_1,1}) (1-q_{l_2,2}) \dots
 (1-q_{l_{s-1},s-1})\cdot q_{l_{s},s}
$$
where ${\cal L}_{s,j}$ is the set of all paths $l=(0,\pm1,*,\dots,*,j\pm 1,j)$ of SRW on $\Z^1$ of length~$s$ ending at point~$j$. Using~\eqref{eq:qbounds} we get that
\begin{align}\label{eq:pjbounds}
\frac 12 p_{j;2/N} \le p_j(m)\le 2 p_{j;1/N}
\end{align}
where $p_{j;\gamma}$ is the corresponding probability for the process which gets killed with a {\em constant} rate $\gamma\in(0,1)$.
This quantity, however, we can compute.
\begin{lemma}\label{lem:pj}
$$
p_{j;\g}=\sqrt{\frac{\g}{2-\g}}\cdot \left[\frac{1-\sqrt{\g(2-\g)}}{1-\g}\right]^{|j|}.
$$
In particular, if $\g=\frac aN$ where $N$ is large and $a=O(1)$,
$$
p_{j;\g}\sim
\sqrt{\frac{a}{2N}}\cdot \left[1-\frac{\sqrt{2a}}{\sqrt{N}} \right]^{|j|}.
$$
\end{lemma}
\begin{proof}
Let $q_i=q_i^{(j)}$ denote the probability that the random walk gets killed at $j$, provided it starts at point $i\in\Z^1$.
We have the following easy recursion:
$$
q_i=
\begin{cases}
(1-\g)\frac{q_{i-1}+q_{i+1}}{2},& \text{if } i\ne j;\\
\g+(1-\g)\frac{q_{i-1}+q_{i+1}}{2},& \text{if } i=j.
\end{cases}
$$
The characteristic equation $\l^2-\frac{2}{1-\g}\l+1=0$
has the roots
$$
\l_{1}=\frac{1+\sqrt{\g(2-\g)}}{1-\g}>1, \quad
\l_{2}=\frac{1-\sqrt{\g(2-\g)}}{1-\g}=\frac{1}{\l_1}<1.
$$
We have different solutions for $i\ge j$ and $i\le j$; moreover, $q_j$ must go to $0$ as $j\to\pm \infty$; this solutions have to be symmetric around $j$. Therefore, we must have
$q_i=C \l_2^{|i-j|}$. Using the recursion at $j=i$ we obtain
$C=\g+(1-\g) C\l_2$ yielding
$C=\frac{\g}{1-(1-\g)\l_2}=\sqrt{\frac{\g}{2-\g}}$. Consequently,
$p_{j;\g}=q_0=\sqrt{\frac{\g}{2-\g}}\cdot \l_2^{|j|}$.
The rest is a simple calculus.
\end{proof}

\begin{lemma}\label{lem:comblower}
$$
\P\left(\xi_N<\frac 18 N^{3/2}\right)=o(1)
$$
for $N$ large.
\end{lemma}
\begin{proof}
As long as $h_j^\pm \ge N/2$, we know from the RHS of~\eqref{eq:pjbounds} and Lemma~\ref{lem:pj} (with $a=1$) that 
$$
p_j(m)\le \frac{2}{\sqrt N}\left[1-\frac{1}{\sqrt {N/2}}\right]^{|j|}.
$$
Therefore, the probability that there will be at least one particle among the first $N^{3/2}$ ones which gets stuck at column $K_j^\pm$ for $|j|\ge N/2$ is smaller than
\begin{align*}
 N^{\frac 32} \sum_{|j|\ge \frac N2} \frac 2{\sqrt N}\left[1-\sqrt{\frac 2N}\right]^{|j|}&\sim
4N \int_{\frac N2}^{\infty} e^{-x\sqrt{\frac 2N }} dx
\\ &
=4N \int_{1}^{\infty} \, e^{-y\sqrt{\frac N2}} dy
=4\sqrt{2N} e^{-\sqrt{\frac N2}}=o(1).
\end{align*}
For the columns with $j$ from $-N/2$ till $+N/2$ this probability is at most $\frac{2}{\sqrt{N}}$. 
We can thus couple our process with independent Bernoulli trials conducted $\frac 18\, N^{3/2}$ times, which has  the average $\frac 14\, N$.
Hence, by large deviation principle (see e.g.~\cite{HOL}), the probability that amongst particles with indices $m=1,2,\dots,\frac 18 N^{3/2}$  
more than $\frac 13\,N$  get stuck at a particular column $j$ is bounded above by
$e^{-\text{const}\cdot \sqrt{N}}$ and hence the probability that at least one $h_j^{\pm}$, $j<|N/2|$, becomes smaller than $N-\frac 13\, N$ is less than $N\cdot e^{-\text{const}\cdot \sqrt{N}}=o(1)$. Consequently, all the $h^{\pm}_j$ indeed remain higher than $N/2$ for the first $\frac 18\, N^{3/2}$ emitted particles. The statement has been proved.
\end{proof}

\begin{lemma}\label{lem:combupper}
$$
\P(\xi_N>20 N^{3/2})=o(1)
$$
for $N$ large.
\end{lemma}
\begin{proof}
Consider the SRW $W_n$ during the first $N$ steps, assuming it is not killed earlier. We have by the reflection principle
\begin{align*}
\P(\max_{n=1,\dots,N} |W_n|>\sqrt{N})
&\le \P(\max_{n\le N} W_n>\sqrt{N})
+\P(\min_{n\le N} W_n<-\sqrt{N})\\
&= 2 \P(\max_{n\le N} W_n>\sqrt{N})
=4\P(W_N>\sqrt{N})=4\left[1-\Phi(1)\right]=0.63\dots
\end{align*}
and hence with probability at least $\approx 0.36\dots$ the walk stays within $[-\sqrt N,+\sqrt N]$ for the first $N$ steps. At the same time the walk is killed at each step with probability at least $1/N$, hence it does not survive until $N+1$ with probability at least $1- (1-1/N)^N\approx 1-e^{-1}=0.63\dots$, that is, the particle gets stuck in one of the columns $K^{\pm}_j$ with $|j|<\sqrt{N}$. 

Consequently, each particle emitted at the origin with probability at least $0.36\cdot 0.63=0.2268$ gets stuck at point inside
$$
A:=\{(x,y)\in\Z^2:\ |y|<N,\ |x|\le \sqrt{N}\}
$$ 
independently of the past. Since $|A|\le 4 N^{3/2}$, and $20\times 0.2268>4$, by the large deviation principle with probability converging to one, $20 N^{3/2}$ points should suffice to fill up $A$ and hence make the origin sticky.
\end{proof}

\section{Aggregation on $\Z^d$, $d\ge 3$}\label{Sec:d-ge-3}
Assume $d\ge 3$ and let $G=[-N,\dots,N]^3\subset \Z^d$  be a cube of the  $d-$dimensional lattice with the sticky border 
${\sf B}=\{{\bf x}\in G:\ |{\bf x}|\ge N-1\}$.
Let $\xi_N$ be the number of particles emitted before the origin $v_0={\bf 0}$ becomes sticky. Trivially, $N\le \xi_N \le (2N)^d$.

In the analogy with Section~\ref{subsecbox2} we will prove the following lower bound for $\xi_N$ (compare this with~\cite{Ke87} for the case $d\ge 3$.)
\begin{thm}\label{thm:Zd}
There exists a $c_d>0$ such that
$\P\left(\xi_N>  c_d N^{d/2} \right)\to 1$ as $N\to\infty$.
\end{thm}
\begin{proof}
Recall that the Green function in dimension $d$
$$
g({\bf 0},{\bf y})\sim \frac{const}{\Vert {\bf y}\Vert^{d-2}}
$$
(see e.g.~\cite{LL}, Theorem 4.3.1) gives the average lifetime number of visits to $\bf y$ of the SRW on $\Z^d$ starting from point ${\bf 0}$. Since the probability of return to ${\bf y}$ of a SRW starting at ${\bf y}$ is a constant smaller than $1$  independent of ${\bf y}$ (due to the transience of the walk on $\Z^d$, $d\ge 3$), conditioned on the first visit to $\bf y$,  the number of visits to $\bf y$ starting from $\bf x$ has a geometric distribution with the same finite mean for all ${\bf y}$; therefore if $X_n$ denotes a SRW on $\Z^d$ then
\bnn\label{eq:Greenf}
\P(X_n={\bf y}\text{ for some }n\ge 1 \| X_0={\bf 0})\le \frac{c}{\Vert {\bf y}\Vert^{d-2}}
\enn
for some constant $c>0$ and sufficiently large ${\bf y}$.

Recall that $B(k)=\{{\bf x}\in G:\ |{\bf x}|\le k\}$, is the ball of radius $k$  distance around the origin and let 
$$
R(k)=B(k+1)\setminus B(k)
$$
be the ``shell'' of radius $k$. Let $\nu_k$ be the index of the first particle to get stuck on $R(k)$.

Denote by $p_k(n)$  the probability of a particle getting stuck in $R(k)$ for the first time, given there are already $n$ sticky points ${\bf y_1},\dots,{\bf y_n}$ in $R(k+1)$. In order for this event to happen, we need at least that a SRW starting from $\bf 0$ hits a point of $R(k)$ adjacent to one of ${\bf y_i}$'s before reaching the boundary $\sf B$; denote these points ${\bf y_1}',\dots, {\bf y_{n'}'}\in R(k)$. Obviously,  $n'\le d n$  due to the fact that ${\bf y_i}'$ must be adjacent to some  ${\bf y_j}$ and at the same time $|{\bf y_i}'|<|{\bf y_j}|$.

It immediately follows from~\eqref{eq:Greenf} that
\begin{align*}
p_k(n) &\le \P \left(
\text{SRW starting at $\bf 0$ ever reaches the set } \{ {\bf y_1}',\dots,{\bf y_{n'}'}\}\right) 
\\
&\le n' \min_{{\bf y}'\in R(k) }\P \left(
\text{SRW starting at $\bf 0$ ever reaches }  {\bf y'}\right)  \le \frac{c\cdot dn}{k^{d-2}}.
\end{align*}

Suppose the particle with index $\nu_{k+1}$ is the first particle to becomes sticky on $R(k+1)$. If the next $(n-1)$ particles do not get stuck at $R(k)$, the number of particles at $R(k+1)$ becomes at most $n$.
Therefore,  the probability that $\nu_k-\nu_{k+1}>n$, which is equivalent to the event that none of the particles with index $\nu_{k+1}+\ell$, $\ell=1,2,\dots,n$ gets stuck at $R(k)$, is at least
$$
\prod_{\ell=1}^n \left(1-\frac{c \ell}{k^{d-2}}\right)
\ge 1-\frac{cdn(n+1)}{2 k^{d-2} }.
$$
Plugging $n=c_1\sqrt{k^{d-2}}$ for a suitable constant $c_1>0$, we get  
$$
\P\left(\nu_k-\nu_{k+1}> c_1 \sqrt{k^{d-2}}
\| \F_{\nu_{k+1}}\right)>\frac 13
$$
where $\F_{\nu_{k+1}}$ is the history of the process up to the time $\nu_{k+1}$.
Consequently, the random variables $(\nu_{k}-\nu_{k+1})$ for $k=N/2,\dots,N$ can be coupled with independent random variables $\eta_k$ taking value  $c_1 \sqrt{(N/2)^{d-2}}=:c_2 N^{d/2-1}$ with probability $1/3$, and $0$ otherwise, such that $\nu_k-\nu_{k+1}\ge \eta_k$.  This in turn yields
\begin{align*}
\P\left(\xi_N\le c_d N^{d/2}\right)
\le 
\P\left(\nu_{N/2}-\nu_N \le  c_d N^{d/2}\right)
\le 
\P\left(\sum_{\ell=N/2}^{N}\eta_\ell \le c_d N^{d/2}\right)
=o(1)
\end{align*}
as long as $c_d<c_2 /6$ by the standard Chernoff bound.
\end{proof}

\section*{Acknowledgment}
S.V.\ would like to thank Yuval Peres for pointing out the equivalence between our model and OK Corral. S.V.\ research is partially supported by Swedish Research Council grant VR~2014-5157.


\begin {thebibliography}{99}

\bibitem{As-Ga-2014}
Asselah, Amine and Gaudilli\`ere, Alexandre.
Lower bounds on fluctuations for internal {DLA},. 
Probab. Theory Related Fields.~158 (2014), no.~1-2, 39--53.

\bibitem{As-Ga-2013}
Asselah, Amine and Gaudilli\`ere, Alexandre.
Sublogarithmic fluctuations for internal {DLA},. 
Ann. Probab.~41 (2013), no.~3A, 1160--1179.

\bibitem{As-Ra-2016}
Asselah, Amine and Rahmani, Houda.
Fluctuations for internal {DLA} on the comb,. 
Ann. Inst. Henri Poincar\'e Probab. Stat.~52 (2016), no.~1, 58--83.


\bibitem{BR}Bhattacharya, R.~N.\ and Ranga Rao, R. Normal approximation and asymptotic expansions.  Wiley Series in Probability and Mathematical Statistics. John Wiley \& Sons, New York-London-Sydney, 1976,

\bibitem{BIL}Billingsley, Patrick.
Probability and measure. 
Second edition. John Wiley \& Sons, Inc., New York, 1986.


\bibitem{BD} Davis, Burgess. Reinforced random walk.  Probab.~Theory Related Fields 84 (1990), no.~2, 203--229.  

\bibitem{DS} Doyle, Peter and Snell, Laurie. Random walks and electric networks.  Carus Mathematical Monographs, 22. Mathematical Association of America, Washington, DC, 1984. 

\bibitem{DUR} Durrett, Rick.
Probability: theory and examples. 
Fourth edition. Cambridge University Press, Cambridge, 2010.

\bibitem{FRE} Freedman, David. Bernard Friedman's urn.  Ann.~Math.~Statist 36, (1965) 956--970. 

\bibitem{HOL} den Hollander, Frank. Large deviations. Fields Institute Monographs, 14. American Mathematical Society, Providence, RI, 2000.


\bibitem{Je-Le-Sh-2013} Jerison, David; Levine, Lionel and Sheffield, Scott. Internal DLA in higher dimensions, Electron. J. Probab.~18 (2013), no.~98, 1083--6489.

\bibitem{Je-Le-Sh-2014} Jerison, David; Levine, Lionel and Sheffield, Scott. Internal DLA and the Gaussian free field,  Duke Math. J.~163 (2014), no.~2, 267--308.

\bibitem{K99} Kingman, J.~F.~C. Martingales in the OK Corral. Bull.\ London Math.\ Soc.\ 31 (1999), no.~5, 601--606. 

\bibitem{KV} Kingman, J.~F.~C.\ and Volkov, S.~E.
Solution to the OK Corral model via decoupling of Friedman's urn. J.~Theoret.~Probab.~16 (2003), no.~1, 267--276. 

\bibitem{LaLy99} Larsen, Michael and Lyons, Russell. Coalescing particles on an Interval.  J.~Theoret.~Probab.~12 (1999), no.~1, 201--205.


\bibitem{LA-Br-Gr-92} Lawler, Gregory; Bramson, Maury and Griffeath, David. Internal diffusion limited aggregation,  Ann. Probab.~20 (1992), no.~4, 2117--2140.

\bibitem{LL} Lawler, Gregory and Limic, Vlada. Random walk: a modern introduction. Cambridge Studies in Advanced Mathematics, 123. Cambridge University Press, Cambridge, 2010.

\bibitem{LE&Pe07} Levine, Lionel; Peres, Yuval. Internal Erosion and the Exponent $\frac 34$. (2007).
\begin{verbatim}
http://www.math.cornell.edu/~levine/erosion.pdf
\end{verbatim}




\bibitem{Ke87} Kesten, Harry. How long are the arms in {DLA}? J. Phys. A.~20 (1987), no.~1, L29--L33.

\bibitem{Ke87_hit} Kesten, Harry. Hitting probabilities of random walks on {${\bf Z}^d$}.  Stochastic Process. Appl.~25 (1987), no.~2, 165--184.

\bibitem{Spit} Spitzer, Frank. Principles of random walk. Second edition. Graduate Texts in Mathematics, Vol.~34. Springer-Verlag, New York-Heidelberg, 1976.


\bibitem{WM} Williams, David and McIlroy, Paul. The OK Corral and the power of the law (a curious Poisson-kernel formula for a parabolic equation). Bull.~London Math.~Soc.~30 (1998), no.~2, 166--170. 


\bibitem{We-Sa-83} Witten, T. A.\ and Sander, L. M. Diffusion-limited aggregation.  Phys. Rev. B (3).~27 (1983), no.~9, 5686--5697.

\end {thebibliography}
\end{document}